\documentclass[11pt]{elsarticle}
\usepackage[margin=1in]{geometry}
\usepackage{lscape}

\usepackage{amsmath,amsfonts,amssymb,amsthm}
\usepackage{graphicx}
\usepackage{xcolor}
\usepackage{comment}
\usepackage{soul}
\usepackage{booktabs}
\usepackage{multirow}

\usepackage{bm} 

\def\br{\br}
\newcommand{\bbx}{\boldsymbol{x}}
\newcommand{\bby}{\boldsymbol{y}}
\newcommand{\bbz}{\boldsymbol{z}}

\newcommand{\gam}{\gamma}

\def\hangone{\vskip0pt\hangindent=28pt\hangafter=1}

\def\hangthree{\vskip0pt\hangindent=48pt\hangafter=1\qquad}

\newtheorem{theo}{Theorem}[section]

\newtheorem{re}[theo]{Remark}

\definecolor{ForestGreen}{rgb}{0.0, 0.5, 0.0}

\begin{document}

\begin{frontmatter}
\author[label1]{Shiping Zhou}
\author[label2]{Yanzhi Zhang
\corref{cor2}}
\cortext[cor2]{Corresponding author: {zhangyanz@umsystem.edu}.
}
\author[label3,label4]{Max Gunzburger}

\address[label1]{Department of Computational Mathematics, Science and Engineering, \\Michigan State University, East Lansing, MI 48824, USA}
\address[label2]{Department of Mathematics and Statistics, Missouri University of Science and Technology, Rolla, MO 65409, USA}
\address[label3]{Department of Scientific Computing, Florida State University, Tallahassee,  FL 32304, USA}
\address[label4]{Oden Institute for Engineering \& Sciences, University of Texas, Austin, TX 78712, USA}

\title{
{
\Large\bf Nonlocal modeling of spatial fractional diffusion \\
with truncated interaction domains and
\\
truncated kernel function singularity   
}
}
\begin{abstract}

\noindent
Parabolic partial differential equations (PDEs) are in ubiquitous, very effective use to model diffusion processes. However, there are many applications (e.g., such as in hydrology, animal foraging, biology, and light diffusion just do name a few) for which results obtained through the use of parabolic PDEs do not agree with observations. In many situations the use of fractional diffusion models {has} been found to be more faithful to that which is observed. Specifically, we replace the Laplacian operator in the PDE by a fractional Laplacian operator ${\mathcal L}$ which is an integral operator for which solutions are sought for on all of space, has an unbounded domain of integration, and for a given point $\bbx$ the integrand contains a kernel function $\phi(\bby-\bbx)$ that is infinite whenever $\bby=\bbx$.

These three features pose impediments not only for the construction of efficient discretization methods but also because all three involve one or more sort of ``infinity''. To overcome these impediments we choose to invoke one or more of the following strategies. 

\hangone(a) We seek solutions only within a chosen bounded domain $\Omega$. 

\hangone(b) For every $\bbx\in\Omega$, we choose a bounded domain of integration such as, e.g., an Euclidean ball {$B_\delta(\bbx)$} having finite radius $\delta$. 

\hangone(c) We truncate the singularity of $\phi(\bby-\bbx)$ by setting, for a given constant $\varepsilon>0$,  $\phi(\bby-\bbx)= \phi(\varepsilon)$ whenever $|\bby-\bbx|\le\varepsilon$.

\noindent We then provide extensive illustrations of the possible combinations chosen from among (a), (b), and (c). We also illustrate, for the various models defined for each of these combinations, their limiting behavior of solutions such as showing that as $\delta\to0$ we recover the PDE model and also showing that in the limit of some other parameters we recover the fractional Laplacian model. 

\end{abstract}

\end{frontmatter}

\section{Introduction}
\label{sec1}

Parabolic partial differential equations (e.g., such as $\partial u/\partial t = \Delta u$ with $\Delta$ denoting the Laplacian operator) are in ubiquitous use for the modeling of diffusion processes; an exemplar of such a process is the transfer of heat from a hot area to a cold area in a material. Associated with parabolic partial differential equations {(PDEs)} are, e.g, the Fourier, Fick, and Ohm laws of diffusion. A long list and discussions about these and other such processes are provided in \cite{Dautray1990} and a rigorous treatment can be found in, e.g., \cite{Dautray1990,Evans2010}. 

However, there are many diffusion processes that cannot be modeled through the use of such laws and therefore through the use of parabolic {PDEs}. Examples of such alternate processes arise in many application areas including  light diffusion \cite{Barthelemy2008}, hydrology \cite{Dentz2004,Kirchner2000}, animal foraging \cite{Kolzsch2015}, and biology \cite{Kusumi2005}, among a myriad of other examples. Collectively, such processes are labelled as  {\em anomalous diffusion} and, nowadays, such models have been successfully used to model a large variety of phenomena that cannot be accurately described by classical {PDEs}.\footnote{As just noted, the term ``anomalous diffusion'' is in ubiquitous use in a myriad of applications. Here, we use that term whenever model solutions do not agree with solutions of PDE diffusion models.}

The invocation of spatial fractional derivatives defines one class of anomalous diffusion models in which a classical {PDE-based} diffusion {model} such as $\displaystyle {\partial u}/{\partial t} = \Delta u$ is replaced by 
\begin{equation}\label{fractfracx}
\qquad \mbox{$\displaystyle \frac{\partial  u}{\partial t} = \big(-(-\Delta)^s \big)u$} 
\end{equation}
\noindent where $-(-\Delta)^s$ denotes the fractional Laplacian operator and $s\in(0,1)$. See, e.g. \cite{DiNezza2012,Metzler2000,Stinga2019} for detailed considerations of models such as that in \eqref{fractfracx}. In Section \ref{sec22} we provide in \eqref{fractionalheat2} and  \eqref{nonLa} two equivalent definitions for the {classical fractional} Laplacian operator $-(-\Delta)^s$. Of note is that the fractional Laplacian operator is an integral operator.

\begin{re}
{\em On a historical note, fractional derivatives were first introduced by Leibniz and then followed by Abel who developed a fractional calculus.   
} \quad $\Box$
\end{re}

Compared to classical {PDE-based} diffusion models, determining analytical solutions {or even numerical solutions} of fractional diffusion equations such as  \eqref{fractfracx} and also understanding and illustrating the properties of those solutions are challenging endeavors. Specifically, two challenges are posed.

\hangone Challenge 1. 
All types of fractional diffusion operators are {\em nonlocal} in nature by which we mean that two points separated by a non-vanishing distance interact with each other.\footnote{In contrast, for PDE models two points interact with each only within infinitesimal neighborhoods needed to define derivatives; see Subsection \ref{sec21}.}  For example, the classical fractional Laplacian operator is an integral operator defined for all $\bbx\in{\mathbb R}^d$ and 
which has a domain of integration which encompasses all $\bby\in{\mathbb R}^d$; see Subsection \ref{sec22}.  

\hangone Challenge 2. The integrand of  fractional Laplacian operators includes a radial kernel function $\phi_s({|\bby-\bbx|})$ that is  {\em singular}, i.e., that is infinite whenever ${\bby}={\bbx}$.

\noindent See, e.g, {\cite{DElia2013,Du2012,Du2013}} for detailed considerations of these challenges.

For Challenge 1, to avoid integration over ${\mathbb R}^d$ and to avoid interactions between two points that are separated by unboundedly large distances, a {\em spatial domain truncation of the kernel function} is introduced; see, e.g., \cite{Acosta2017,Burkovska2019,DElia2013}. Specifically, the kernel function is set to zero whenever $|\bby-\bbx|>\delta>0$, i.e., for any $\bbx$, the kernel function is nonzero only for $\bby\in B_\delta(\bbx)$, where $B_\delta(\bbx)$ denotes the Euclidean ball centered at $\bbx$ having radius $\delta$; here $\delta$ is a given bounded parameter that is often referred to as the {\em horizon} or the {\em interaction radius} and $B_\delta(\bbx)$ is referred to as the  {\em interaction neighborhood}, among other monikers. Thus, because we now have that points $\bbx\in\Omega$ only interact with points $\bby\in B_\delta(\bbx)$, we have that, for a given domain $\Omega$, {\em the domain of integration reduces to the finite domain $\Omega\cup\Omega_{\mathcal I}$, where {$\Omega_{\mathcal I} = \{ \bby\in B_\delta(\bbx) \cap ( {\mathbb R}^d \setminus \Omega)$ for all $\bbx\in\overline\Omega$\}}.} 

With respect to Challenge 2, one can avoid the singularity by a {\em kernel function truncation}, i.e., the singularity in the kernel function $\phi_s(|\bby-\bbx|)$ occurring at $\bby=\bbx$  {is removed  by setting, for a given constant $\varepsilon$, the kernel function to be a suitable bounded function whenever $|\bby-\bbx|\le\varepsilon$.}

In this paper we focus on {\em seeking solutions of fractional diffusion equations in bounded domains $\Omega\subset {\mathbb R}^d$} but not on all of ${\mathbb R}^d.$

\vskip5pt
\begin{re}
{\em
For the sake of clarity, we henceforth keep the notation simple by {\em confining ourselves to one spatial dimension}; all that is written for this case easily extends to multiple spatial dimensions. In addition, we confine ourselves to rather simple linear diffusion processes; again, extensions to other types of fractional models are straightforward.
}\quad$\Box$
\end{re}

The paper is organized as follows. We first review, in Section \ref{sec2}, parabolic PDE models for which solutions are sought for on a bounded domain $\Omega=(-L,L)$, fractional diffusion models for which solutions are sought for on all of ${\mathbb R}^d$, and fractional diffusion models for which solutions are sought for on a bounded domain $\Omega=(-L,L)$.

In Section \ref{sec3} we first define approximations of the kernel function $\phi_s(|y-x|)$ constructed by truncating its spatial support or by truncating its singularity or both.\footnote{To be clear, even though we confine ourselves to seeking solutions at points $x$ in the bounded domain $(-L,L)$, in the fractional system \eqref{fractionalheat5} the support of the kernel function is infinite. Hence, in  Section \ref{sub31}, we  truncate the support of the kernel function to the finite interval $[-\delta,\delta]$.} We then consider two fractional heat equations posed on points within a bounded spatial {domain} and which have kernel functions that have truncated spatial support and have either singular kernel functions or are truncated bounded kernel functions. 

In Section \ref{sec4} we consider four fractional diffusion systems for all of which solutions are sought in a bounded domain $\Omega=(-L,L)$; {naturally, such a restriction is of most interest in practice.} Specifically, the four systems feature kernel functions which have bounded spatial support $[-\delta,\delta]$ for two values of $\delta$ and for which the singularity is or is not truncated. In so doing, we systematically study the effects that have on solutions due to the use of truncated support of the spatial domain and due to the use of truncated singularities of kernel functions. 

{

For convenience, we introduce the following notations that are used in the paper. Specifically, we have the input parameters, functions, and input data given as follows:
\begin{equation}\label{list}
\begin{tabular}{llllr}
\quad
&  
& spatial interval for solutions $u(x)$
& $x\in [-L,L] \subseteq{\mathbb R}$\quad for\,\, $L>0$ 
\\[.5ex]
\quad
& 
& temporal interval for $u(x)$
& $ t\in[0,T]\subseteq  {\mathbb R}\backslash{\mathbb R}^-$
\\[.5ex]
\quad
& 
& initial condition data 
& $u_0(x)$\quad for \,\, $ x \in [-L,L]$ 
\\[.5ex]
\quad
& 
& volume constraint data 
& $g_-(x,t)$\quad for \,\,$x \in (-\infty,-L]$\,\, and \,\, $t > 0$
\\
\quad
&  
&  
& $g_+(x,t)$\quad for \,\,$x \in [L, \infty)$\,\, and \,\, $t > 0$ \\[.5ex]
\quad
& 
&singular kernel function
& for any\,\, $z \in {\mathbb R}$ \,\,\,and\,\,\,$s\in(0,1)$ 
&  
\\  
&&&
\quad $\displaystyle
\phi_s(z) =  \frac{C_{s}}{|z|^{1+2s}}$ \,\,\,with\,\,\,
$\displaystyle C_s = \frac{ 4^s s \Gamma(s+1/2) }{ \pi^{1/2} \Gamma(1-s)   }$
\\[.5ex]
\quad
& 
& {spatial support truncation radius} 
& $\delta \in (0,\infty)$ 
\\[.2ex]
\quad
& 
& {singularity truncation radius} 
& $\varepsilon >0$
\\[.2ex]
\quad
& 
& {spatial grid size} 
& $h$
\\[.2ex]
\quad
& 
& {temporal time step} 
& $\Delta t$
\end{tabular}
\end{equation}
\noindent
Here $\Gamma(\cdot)$ denotes the Gamma function and $L=\infty$  and $T=\infty$ are allowable.\footnote{In ${\mathbb R}^d$ we would have that the kernel function is given by {$\phi_{d,s}({\bbz})=C_{d,s}\big/{|{\bbz}|}^{d+2s}$} \\with {$C_{d,s}=4^s s \Gamma(s+d/2) \big/\big(\pi^{d/2} \Gamma(1-s)\big)$.}} Note that the kernel function $\phi_s(z)$ is infinite at $z=0$ and is a monotone decreasing function for increasing values of $|z|$. 

\begin{re}\label{rerere}
{\em
As we repeatedly observe in the rest of the paper, the value of $s$ has a profound effect on the spatial and temporal evolution of the solution $u(x,t)$.
}\qquad$\Box$
\end{re}

\subsection{\bf About the discretization of the systems considered}\label{sec11}

The challenges discussed above are concerned with continuous fractional diffusion models such as \eqref{fractfracx}. Certainly, because exact solutions of continuous models are seldom available, accurate numerical simulations are, of course, necessary to obtain useful information. Here, all plots appearing in the paper were generated using accurate computational approximations. 
For all cases considered, temporal discretization is effected via the Crank--Nicolson method. Unless otherwise noted, a uniform time-step size $\Delta t=0.0001$ is used. 

Two types of spatial discretizations are used. 

\hangone-- For the fractional systems considered in Subsections \ref{sec23}, \ref{sys42}, and \ref{sys43} the spectral method of \cite{Zhou2024} is used; this method handles the singularity of the kernel function analytically without any approximations. For the system in  Subsection \ref{sec23} the spatial support of the kernel function is unbounded whereas for the systems in Subsections \ref{sys42} and \ref{sys43} the spatial support of the kernel function is bounded.

\hangone-- For the spatial discretization of the PDE system considered in Subsection \ref{sec21} a standard finite difference method is used. Also, for the systems considered in Subsections \ref{sys44} and \ref{sys45} the kernel function support domains and singularities are truncated and the finite difference method for nonlocal problems developed in  \cite{Duo2019-Comp, Duo2018, Duo2019-TFL,DElia2020,DElia2013} 
is used. Unless otherwise noted, all three systems employ a uniform grid-size $h=0.0025$ and solutions are sought for in the one-dimensional spatial domain $\Omega = (-L, L)$.

With respect to the initial condition, we would ideally like to impose the Dirac delta function at $x=0$ as the initial condition for all the systems we consider. Such an initial condition would result in transparent behaviors of the temporal evolution of the solution $u(x,t)$. However, of course, for the computational simulations one has to avoid the singularity of the Dirac delta function. Thus for all our simulations as an initial condition we employ a ``smoothed'' Dirac delta function which is large at $x=0$ and then decays quickly as $x$ increases. Specifically we set
\begin{equation}\label{initcond}
u_0(x)=u(x, 0) = \frac{1}{\sqrt{4\pi \widetilde t}}\exp\Big(-\frac{x^2}{4{\widetilde t}}\Big)
\quad \mbox{with} \,\,\, {\widetilde t} =  0.0001.
\end{equation}
Table \ref{tab:C11111} shows that the initial condition defined in \eqref{initcond} decays very quickly as $x$ increases.
\begin{table}[htb!]
\begin{center}
\begin{tabular}{|c||ccccc|}
\hline
$x =$& 0 & $0.0025=h$ & $0.01=4h$ & $0.05=20h$& $0.1=40h$  \\
\hline
$u_0(x)=u(x,0)\approx$ &28.2095 &27.7721 &21.9696 & 5.4457e-2 &3.9177e-10  \\
\hline
\end{tabular}
\caption{The initial condition $u_0(x)=u(x,0)$ in \eqref{initcond} for five values of $x$ with $h=0.0025$.}
\label{tab:C11111}
\end{center}
\end{table}

\section{Partial differential equation and fractional equation diffusion systems}\label{sec2}

In this section we introduce three models for diffusion processes. In the subsequent sections we provide alternative approximations of some of those models that can be realized computationally.

\subsection{\bf Classical partial differential equation system for heat conduction}\label{sec21}

A classical heat equation system governing the spatial and temporal behavior of a function $u^{\it pde}(x,t)$ is given by\footnote{Of course, general classical heat equations also feature a given positive constitutive function $\gam(x,t)$ and a given right-hand side function $f(x,t)$, i.e., we have $u_t - (\gam(x,t) u_{x})_x = f(x,t)$; see, e.g., \cite{Evans2010}, for details. We do not include these functions in our deliberations because they are not germane to our goals and also to simplify the exposition.}
\begin{equation}\label{classicalheat}
\left\{
\begin{aligned}
\frac{\partial u^{\it pde}}{\partial t}(x,t) = \frac{\partial^2 u^{\it pde}}{\partial x^2}(x,t)  &\qquad\mbox{for} \,\,\, x\in (-L,L)   \,\,\,\mbox{and}\,\,\,\,  t\in(0,T] \\
u^{\it pde}(-L,t) = g_-(-L,t)   &\qquad\mbox{for}\,\,\,\, t\in(0,T] \quad \\
u^{\it pde}(L,t) = g_+(L,t)  &\qquad\mbox{for}\,\,\,\, t\in(0,T] \quad\\ 
u^{\it pde}(x,0) = u_0(x) 
 &\qquad\mbox{for} \,\,\, x\in [-L,L] \\
\end{aligned}
\right.
\end{equation}
{where the spatial interval, temporal interval, and initial conditions are defined in \eqref{list} and we are given the boundary condition input data $g_-(-L,t)$ and $g_+(L,t)$ for $t\in(0,T]$ and the initial data $u_0(x)$ for $x\in[-L,L]$. The choice $L=\infty$ is allowable in which case the second and third lines of \eqref{classicalheat} are omitted.

It is well known that if the boundary conditions are homogeneous, i.e., $g_-(-L,t)=0$ and $g_+(L,t)=0$, then, after short a temporal interval, the solution $u^{\it pde}(x,t)$ of \eqref{classicalheat} decays pointwise in $x$ and also in norm; specifically we have that
\begin{equation}\label{decayclassicalheat}
|u^{\it pde}(x,t)| \propto \frac{1}{t^{1/2}}  \qquad  \mbox{and}\qquad \|u^{\it pde}(\cdot,t)\|_{L^2{([-L,L])}}  \propto \frac{1}{t^{1/4}}.
\end{equation}
These decay rates also holds for the case $L=\infty$ provided that the initial condition data $u_0(x)$ satisfies certain boundedness conditions.

\begin{figure}[htb!]
\centerline{
\raisebox{6mm}{(a)}
\includegraphics[width=0.42\linewidth, height = 5.6cm]{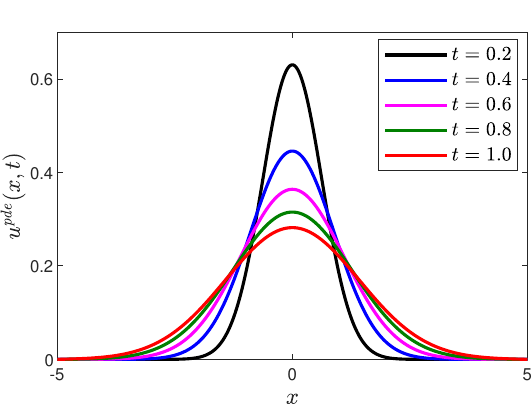}\hspace{2mm}
\includegraphics[width = 0.42\linewidth, height = 5.3cm]{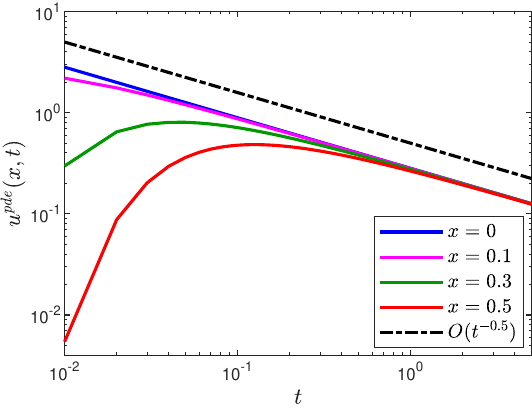}\raisebox{6mm}{\quad(b)}
}
\vspace{-2mm}
\caption{System \eqref{classicalheat}.
{(a)} Plots of $u^{pde}(x,t)$ at five time instances; the initial condition is given in \eqref{initcond}. The decay of the peak and the spread of the solution as time increases are patently clear.
{(b)} Plots of $u^{pde}(x,t)$ for four values of $x$ that illustrate the first formula in \eqref{decayclassicalheat} which, as $t$ increases, conform to the slope of the black dashed line that has the slope of that formula.
}\label{temp-pde}
\end{figure}

Figure \ref{temp-pde}(a) illustrates the decay and spread of the solution $u^{\it pde}(x,t)$ as time increases. Figure \ref{temp-pde}(b) shows that} the temporal slope of the value of $u^{\it pde}(x,t)$ decays in time according to the first formula in (\ref{decayclassicalheat}). For $x\ne0$, the slopes eventually ``catch up'' to the theoretical slope as is illustrated for three values of $x\ne0$.\footnote{Note that many of the figures in the paper, including Figure \ref{temp-pde}{(b)}, are log-log plots.}

\subsection{\bf Classical fractional heat equation system}\label{sec22}

We next recall the classical fractional heat equation system given by\footnote{Details that are germane to this subsection and Subsection \ref{sec23}  can be found in, e.g., \cite{DiNezza2012,Metzler2000,Stinga2019,Vazquez2018}.}
\begin{equation}\label{fractionalheat1}
\mbox{for $s\in(0,1)$} \quad \left\{
\begin{aligned}
u_t^{s}(x, t) = {\mathcal L}_{s} u^{s}(x, t)&\qquad\mbox{for} \ x\in  {\mathbb R} \,\,\,\mbox{and}\,\,\, t\in  (0,T] \quad \\
u^{s}(x,0) = u_0(x) &\qquad\mbox{for} \ x\in{\mathbb R}\\
\end{aligned}
\right.
\end{equation}
where the {\em integral form} of the fractional Laplacian ${\mathcal L}_s$ operating on a function $u^s(x, t)$ is, for $x\in{\mathbb R}$ and {$t > 0$}, defined as
\begin{equation}\label{fractionalheat2} 
 {\mathcal L}_{s} u^{s}(x,t) = {\big(- (-\Delta)^s\big)} u^{s}(x,t) 
 =  \mbox{\em P.V.}  \int_{{\mathbb R}} \Big(u^{s}(y,t)-u^{s}(x,t)\Big)\phi_s(|y-x|) dy 
\quad \mbox{for} \,\,\,s \in (0,1)
\end{equation}
{where {\em P.V.}  denotes the Cauchy principal value and $\phi_s(\cdot)$ is defined in \eqref{list}. 
Equivalently, the {\em spectral form} of the fractional Laplacian ${\mathcal L}_s$ operating on a function $u^s(x, t)$ is, for {$t > 0$}, defined as 
\begin{equation}\label{nonLa}
\begin{aligned}
	{\mathcal L}_{s} u^{s}(x,t)
	&= {\big(- (-\Delta)^s\big)} u^{s}(x,t)
	=  -{\mathcal F}^{-1}{\Big(}|\xi|^{2s} \big({\mathcal F} \,{u^s}(x,t)\big)\Big)
	 \qquad \mbox{for \ $s\in(0, 1)$}
\end{aligned}
\end{equation}
with
$$
{\mathcal F} u(x,t) = \frac{1}{(2\pi)^{1/2}} \int_{\mathbb R} e^{-i \xi x} {u(x,t)} d\xi
$$
denoting the Fourier transform.  
{This definition in \eqref{nonLa} is actually valid for any $s > 0$; however here we are only interested in $s \in (0, 1)$.} 
}

It is clear that the fractional Laplacian is a  {\em nonlocal operator} in the sense defined in {\cite{DElia2020,DElia2013,Du2012,Du2013,Metzler2000}, i.e., points $x\in{\mathbb R}$ interact with all points $y\in{\mathbb R}$.
Also, a discussion of the physics and applications of the fractional Laplacian in the anomalous diffusion setting can be found in \cite{Metzler2000}.

{It is known (see, e.g., \cite{DiNezza2012,Metzler2000}) that for a regular function $u^s(x,t)$ the fractional Laplacian operating on  $u^{s}(x,t)$  satisfies  }   
\begin{equation}\label{fractionalheat3}
\left\{\begin{aligned}
&\mbox{as}\,\,\, s\to 0^+ & {\mathcal L}_{s} u^{s}(x,t) \rightarrow  - u^{s}(x,t)   \\
&\mbox{as}\,\,\, s\to 1^-  &{\mathcal L}_{s} u^{s}(x,t) \rightarrow \Delta u^{s}(x,t) 
\end{aligned}\right.
\end{equation}
where $\Delta$ denotes the classical Laplace partial differential operator.

\vskip5pt

\begin{re}
{{\em The notation $0^+$ and $1^-$ which is often used in this paper respectively denote generic values that are greater than but very nearly zero and generic values that are smaller than but very nearly one.}}\quad$\Box$
 \end{re}

It is also known that, for a given initial condition $u_0(x)$, the solution of \eqref{fractionalheat1} is given by
\begin{equation}\label{fractionalheat4444}
  u^{s}(x,t) = \int_{\mathbb R} \psi_s(x-y,t) u_0(y) dy
\end{equation}
where $\psi_s(x,t)$ denotes the inverse Fourier transform corresponding to the function $f_s(\xi)=e^{-|\xi|^{2s}t}$ so that determining a solution of the system \eqref{fractionalheat1} reduces to the evaluation of the integral \eqref{fractionalheat4444}; see \cite{Vazquez2018}.  For $s=1$, we have the PDE case corresponding to the first kernel function $\psi_{pde}(x,t)$ in \eqref{fractionalheat5555}; that kernel function is simply the Gaussian kernel function for the classical heat equation. 

For $s\ne 1$ the function $\psi_s(x, t)$ can be defined analytically in terms of the Fox $H$-function \cite{Mainardi2005}. Hence, using high-accuracy quadrature rules, very accurate approximations of the exact solution of the system \eqref{fractionalheat1} can be obtained. Such accurate solutions can then be used to determine the error incurred by the solution of discrete approximations of \eqref{fractionalheat1} and \eqref{fractionalheat2}. For the special case of $s=1/2$, $\psi_{1/2}(x, t)$ is given by the right formula in \eqref{fractionalheat5555}:
\begin{equation}\label{fractionalheat5555}
\psi_{pde}(x,t)= (4\pi t)^{- 1/2} e^{-x^2/4t} \,\,\,\mbox{for $s=1$}
\qquad\mbox{and}\qquad
\psi_{1/2}(x,t)= \frac{t}{ \pi(t^2 + x^2)} \,\,\,\mbox{for $s=1/2$.}
\end{equation}

\subsection{\bf Classical fractional heat equation system posed on bounded domains}\label{sec23}

For the classical fractional diffusion system given in \eqref{fractionalheat1} not only does the operator ${\mathcal L}_{s}$ act on a function $u^{s}(x,t)$ {\em for all $x\in\mathbb R$} but it also {requires an integration over all $y\in\mathbb R$.}\footnote{As a result, there is no need to define volume constraints.} 
Here, we consider the case in which an operator ${\mathcal L}_{s,L}$ acts on a function $u^{s,L}(x,t)$ only over the bounded domain $x\in(-L,L)$, i.e, we have, for $s \in (0,1)$ and $t>0$, 
\begin{equation}\label{fractionalheat4}
    {\mathcal L}_{s,L} u^{s,L}(x,t) = {P.V.}\int_{{\mathbb R}} \Big(u^{s,L}(y,t)-u^{s,L}(x,t)\Big)\phi_s(|y-x|)dy 
     \quad\mbox{for}\,\,\,x\in(-L,L).
\end{equation}
Then, the classical fractional diffusion system on bounded domains takes the form
\begin{equation}\label{fractionalheat5}
\left\{
\begin{aligned}
u_t^{s,L}(x, t) = {\mathcal L}_{s,L} u^{s,L}(x, t)&\qquad\mbox{for} \ x\in  (-L,L) \,\,\,\mbox{and}\,\,\, t\in (0,T] 
\\
u^{s,L}(x,t) = g_-(x,t) &\qquad\mbox{for} \,\,\,x\in (-\infty,-L] \,\,\,\mbox{and}\,\,\, t\in (0,T] 
\\
u^{s,L}(x,t) = g_+(x,t) &\qquad\mbox{for} \,\,\,x\in [L,\infty) \,\,\,\mbox{and}\,\,\, t\in (0,T]      
\\
u^{s,L}(x,0) = u_0(x) &\qquad\mbox{for} \,\,\,x\in[-L,L]
\end{aligned}
\right.
\end{equation}
where now it is necessary to impose volume constraints over $x\in {\mathbb R}\backslash(-L,L)$ and the initial condition is imposed only for $x\in[-L,L]$. In \eqref{fractionalheat4} the kernel function has support $(-\infty,\infty)$ and is also singular.  

Note that {unlike the systems considered in the subsequent sections, in this subsection neither the infinite integration domain nor the singularity of the kernel function is truncated.} 

For the fractional  heat system \eqref{fractionalheat5}, {if the volume {constraints} are homogeneous, i.e., $g_-(x,t)=0$ and $g_+(x,t)=0$} {for $x \in {\mathbb R}\backslash(-L, L)$}, our numerical explorations indicate that the temporal decay rates are given by 
\begin{equation}\label{us0rate} 
|u^{s,L}(0, t)| \propto t^{-\frac{1}{2s}} \quad\mbox{and}\quad \|u^{s,L}(\cdot, t)\|_{\ell^2} \propto t^{-\frac{1}{4s}} \quad \ \mbox{for} \ \, s \in (0, 1).
\end{equation}
These results are consistent with the analytical decay rate given in \eqref{fractionalheat5555} for $s = 1/2$ and are also consistent with the limit $u^{s,L}(x, t) \to  u^{\it pde}(x, t)$ as $s\to1^-$. 

We explore, through computational simulations, the effects that the choices made for $s$ and $L$ have on the temporal evolution of the solutions $u^{s,L}(x,t)$. We divide our numerical explorations into two sets because the solution evolutions in ``earlier'' times (say $t\le 1$) differ from those in ``later'' times (say for $t> 1$).

In Table \ref{tab:tbleone} we provide the peak values of $u^{s, L}(0, t)$ for six time instances {$t\le 1$} and for five different $s$ values and also of the peak value of $u^{pde}(0, t)$. For every $s$ and for the PDE case, the peak values decrease monotonically as the time instances increase. Also, for the same value of $t$, the peak values decrease monotonically as $s$ increases and the PDE solution has an even smaller peak value.

\begin{table}[h!]
\begin{center}
\begin{tabular}{|l||cccccc|}
\hline
$t =\,\longrightarrow$& 0 & 0.2 & 0.4 & 0.6  &  0.8 & 1   \\
\hline
$s = 0.001$  &28.2095  &23.0624 &18.8544  &15.4143 &12.6018  &10.3025 \\
$ s = 0.25$  & 28.2095 &8.3645  &3.3376 & 1.6658   &0.9710   &0.6292 \\
$ s = 0.50$  & 28.2095 &1.5837  &0.7947 & 0.5301   &0.3976   & 0.3179 \\
$ s = 0.75$  & 28.2095 &0.8397  &0.5292 & 0.4039   &0.3334   & 0.2873 \\
$ s = 0.90$  & 28.2095 &0.6919  &0.4709 & 0.3759   &0.3204   & 0.2830 \\
\hline
\,\quad{\it pde} & 28.2095 &0.6306 &0.4460 & 0.3642 &0.3154   & 0.2821  \\
\hline
\end{tabular}
\caption{For \(L=5\) and six time instances \(t\le 1\), comparison of the peak values \(u^{s,L}(0,t)\) of the fractional diffusion system \eqref{fractionalheat5} for five different \(s\) values and for the peak values $u^{pde}(0,t)$ for the classical diffusion system  \eqref{classicalheat}.
}\label{tab:tbleone}
\end{center}
\end{table}

\begin{figure}[htb!]
\centerline{
   \raisebox{6mm}{(a)}
   \includegraphics[width=0.415\linewidth]{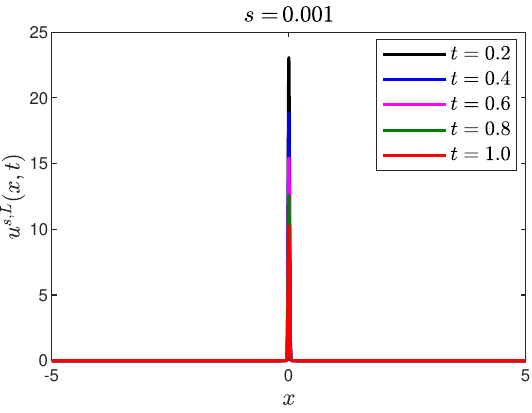}\hspace{5mm}
   \includegraphics[width=0.415\linewidth]{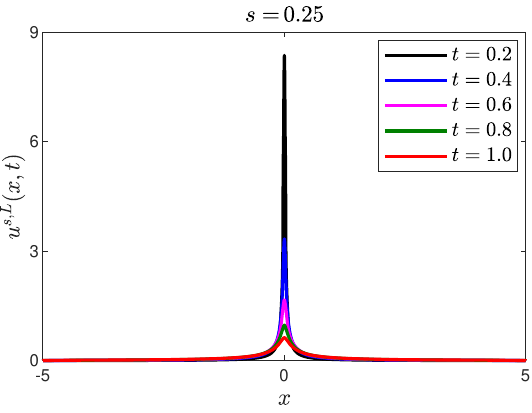}
   \raisebox{6mm}{\quad(b)}
}\vspace{-1mm}
\centerline{
    \raisebox{6mm}{(c)}
    \includegraphics[width=0.415\linewidth]{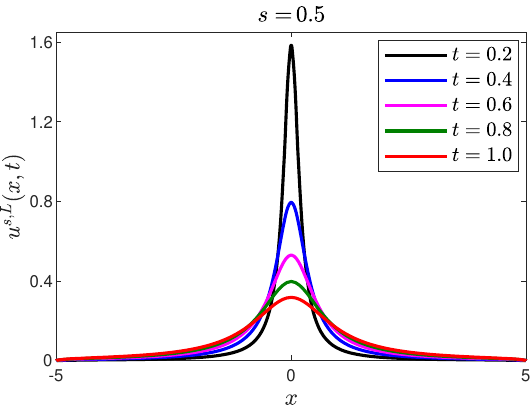}\hspace{5mm}
    \includegraphics[width=0.415\linewidth]{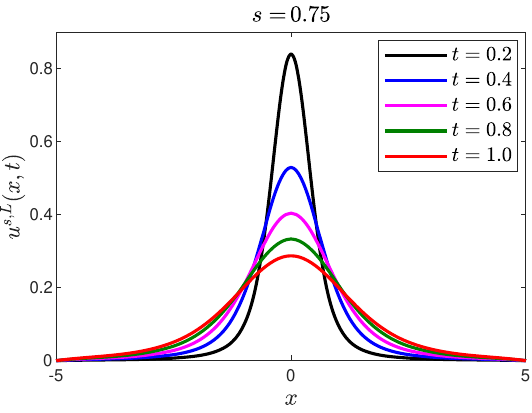}\raisebox{6mm}{\quad (d)}
    }\vspace{-1mm}
    \centerline{
    \raisebox{6mm}{(e)}\includegraphics[width=0.415\linewidth]{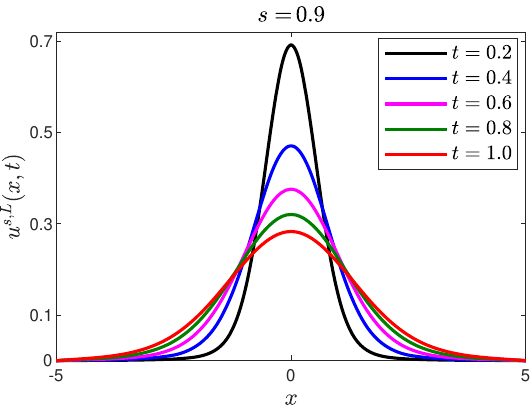}\hspace{5mm}
    \includegraphics[width=0.415\linewidth]{figures-trunc/classical_heat_sol.pdf}\raisebox{6mm}{\quad (f)}
    }\vspace{-3mm}
    \caption{{System \eqref{fractionalheat5}. For $L=5$ and for five values of $s$, the time history of the solution $u^{s,L}(x,t)$ of the fractional heat system \eqref{fractionalheat5} for $t \le 1$. The last plot does the same for the  solution $u^{\it pde}(x,t)$ of the PDE system \eqref{classicalheat}.}}  
\label{fig:fgleone}
\end{figure}

In Figure \ref{fig:fgleone}, for $L=5$ and for five values of $s$, {we present} the time history of the solution $u^{s,L}(x,t)$ of the fractional heat system \eqref{fractionalheat5} { for  $t \le 1$}. The last plot does the same for the  solution $u^{\it pde}(x,t)$ of the PDE system \eqref{classicalheat}. Inspection of the peak values of the plots in Figure \ref{tab:tbleone} are, of course, in agreement with  Table \ref{tab:tbleone}. 
If we examine Figure \ref{fig:fgleone} at say $x=2.5$, we observe that the values of the solutions are now monotonically {\em increasing} as $t$ increases. It is clear that the order of the values has reversed from those at $x=0$. Also, there are points somewhere between $x=0$ and $x=2.5$ where the value of the solutions are not at all monotone.

\begin{figure}[htb!]
    \centerline{
    \raisebox{6mm}{(a)}
    \includegraphics[width=0.41\linewidth]{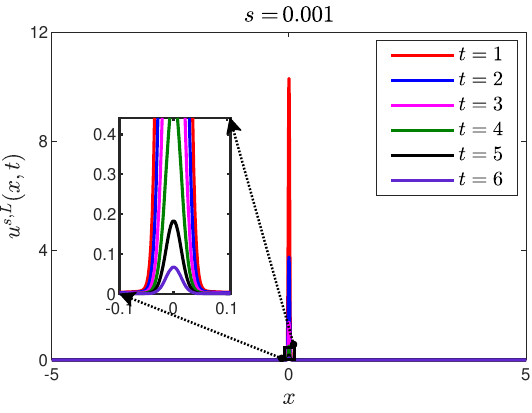}\hspace{5mm}
    \includegraphics[width=0.415\linewidth]{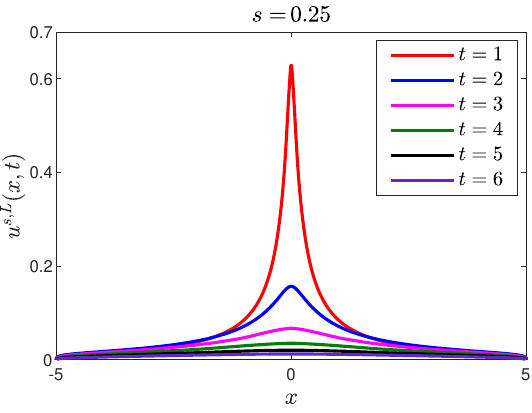}
    \raisebox{6mm}{\quad(b)}}\vspace{-1mm}
    \centerline{
   \raisebox{6mm}{(c)}
   \includegraphics[width=0.415\linewidth]{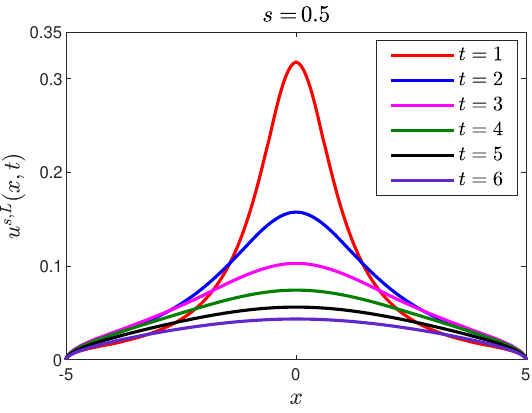}\hspace{5mm}
   \includegraphics[width=0.415\linewidth]{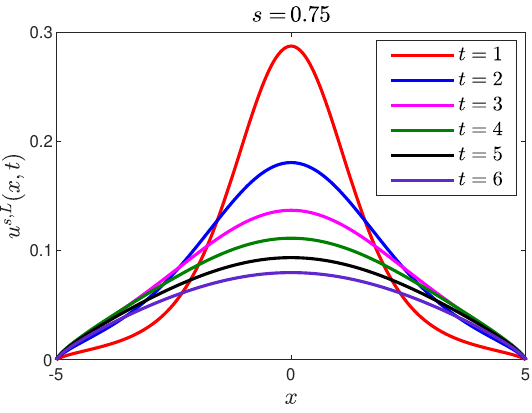}
   \raisebox{6mm}{\quad(d)}
    }\vspace{-1mm}
    \centerline{
   \raisebox{6mm}{(e)}
   \includegraphics[width=0.415\linewidth]{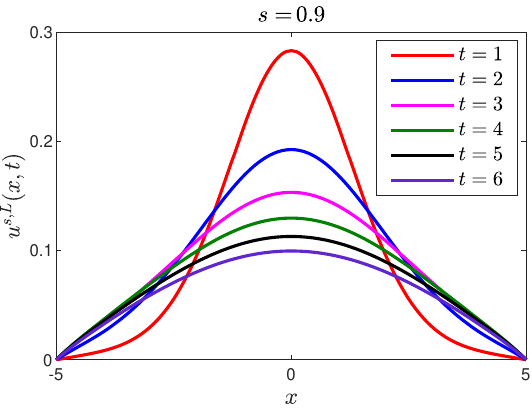}\hspace{5mm}
   \includegraphics[width=0.415\linewidth]{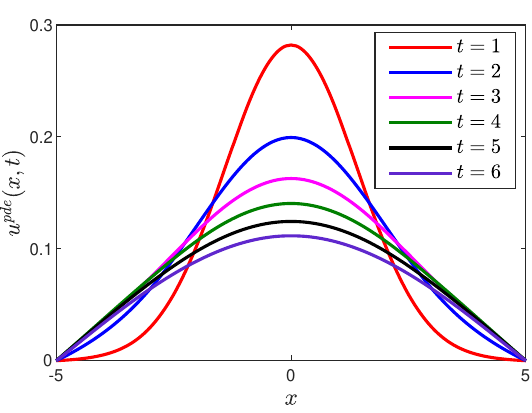}
   \raisebox{6mm}{\quad(f)}
    }\vspace{-3mm}
    \caption{System \eqref{fractionalheat5}. For $L=5$ and for {five} values of $s$, the time history of the solution $u^{s,L}(x,t)$ of the fractional heat system \eqref{fractionalheat5} {for times $t \ge 1$}. The last plot does the same for the  solution $u^{\it pde}(x,t)$ of the PDE system \eqref{classicalheat}.}  
\label{fig:fggeone}
\end{figure} 
In Figure \ref{fig:fggeone} we present the solutions $u^{s,L}(x,t)$ and $u^{\it pde}(x,t)$ for $L = 5$ and for larger times $t \ge 1$. We see that the transitions between the two types of ordering occur at larger values of $x$ as $s$ increases.
\begin{figure}[htb!]
    \centerline{
    \raisebox{6mm}{(a)}
    \includegraphics[width=0.415\linewidth]{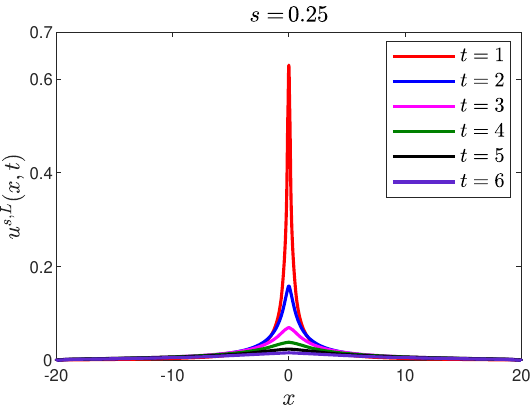}\hspace{5mm}
    \includegraphics[width=0.415\linewidth]{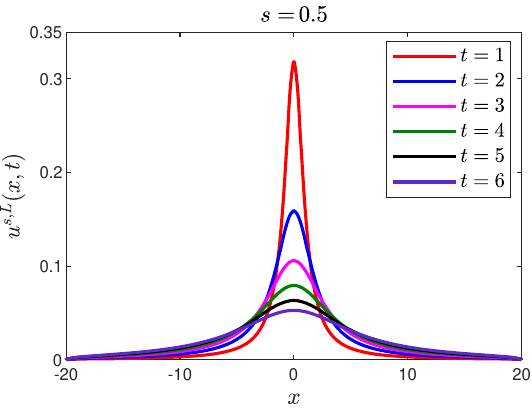}\raisebox{6mm}{\quad(b)}}
    \vspace{-1mm}
    \centerline{
   \raisebox{6mm}{(c)}
   \includegraphics[width=0.415\linewidth]{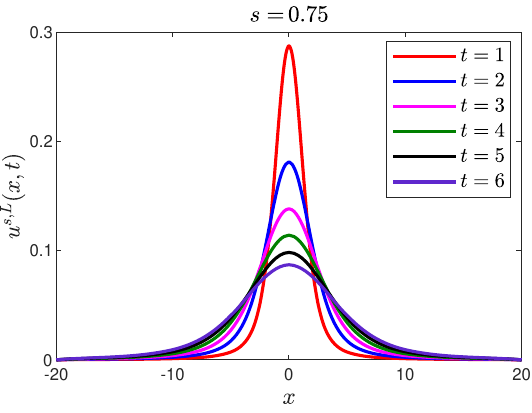}\quad
   \includegraphics[width=0.415\linewidth]{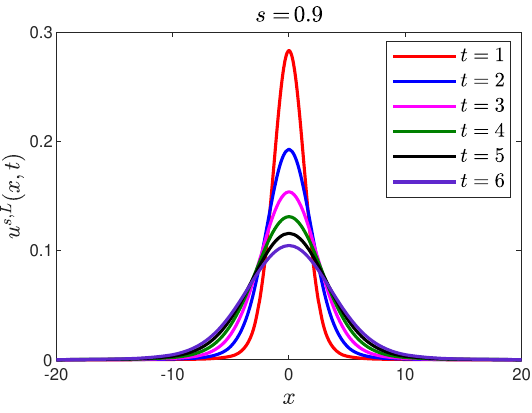}\raisebox{6mm}{\quad(d)}
    }\vspace{-1mm}
    \centerline{\raisebox{6mm}{(e)}
    \includegraphics[width=0.415\linewidth]{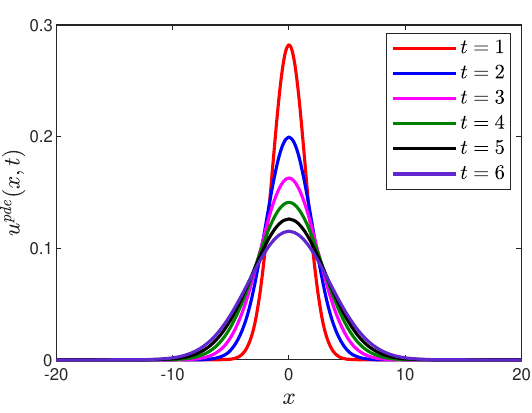}
    }\vspace{-3mm}
    \caption{System \eqref{fractionalheat5}. For $L=20$ and for four values of $s$, the time history of the solution $u^{s,L}(x,t)$ of the fractional heat system \eqref{fractionalheat5} for time $t \ge 1$. The last plot does the same for the  solution $u^{\it pde}(x,t)$ of the PDE system \eqref{classicalheat}.}
    \label{fig:fggeoneL20} 
\end{figure}

We have studied the effects that the values of $s$ and $t$ have on the solutions of the system \eqref{fractionalheat5}. We now consider the effects that $L$, i.e., the length $2L$ the domain of $\Omega$, has on solutions. In Figure~\ref{fig:fggeoneL20}, we present the time evolution of $u^{s,L}(x,t)$ and $u^{\it pde}(x,t)$ for $L=20$ and $t\ge 1$, in comparison to those in Figure \ref{fig:fggeone} for $L = 5$. Clearly the two solutions {in Figures \ref{fig:fggeone} and \ref{fig:fggeoneL20}} are very different. Of course, for $L=5$ all solutions vanish for $x>5$ and $x<-5$ because of the given volume constraints and for $L=20$, the solutions at $x=5$ and $x=-5$ and beyond are certainly appreciable. What is interesting is the interplay between $s$ and $L$. For smallish values of $s$, e.g., $s=0.5$, the solutions at $x=\pm 5$ are much more appreciable than they are for, say, $s=0.9$. Furthermore, for sufficiently large $L$, the fractional heat equation \eqref{fractionalheat5} provides a good approximation to the unbounded system \eqref{fractionalheat1}. However, the longer the simulation time window $(0, T]$, the larger the domain $L$ is needed to suppress boundary effects. 

In Figure~\ref{fig:peak}, we compare the time evolution of the peak value $u^{s,L}(0, t)$ for four choices of $s$ and for  $u^{pde}(0,t)$, with $L = 20$ and $t \in [0, 6]$. For each $s$, the peak decays monotonically in time. Consistent with Table~\ref{tab:tbleone}, at early times (such as $t \le 1$) the peak value decreases monotonically with $s$; in particular, $u^{pde}(0, t)$ has the smallest value. However, for sufficiently large $t$ (such as $t \ge 4$) the ordering reverses: the peak increases with $s$, and $u^{pde}(0, t)$ becomes the largest. This crossover reflects the competition between local smoothing (which strengthens as $s$ increases) and long-range transport (more significant for smaller $s$). Over longer horizons, nonlocal transport for smaller $s$ depletes the center more rapidly, yielding a lower peak.

\begin{figure}[htb!]
    \centering \includegraphics[width=0.43\linewidth]{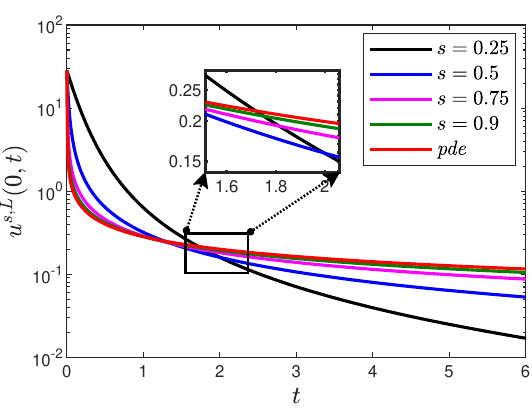}
    \vspace{-3mm}
    \caption{For $L = 20$, the time evolution of the peak $u^{s, L}(0, t)$ of the fractional heat system \eqref{fractionalheat5} and $u^{pde}(0, t)$ of the PDE system \eqref{classicalheat}.}
    \label{fig:peak}
\end{figure}

So far we have concentrated mostly on, for a given $s$, the time histories of solutions. Here, we briefly turn the tables to further consider, for a given $t$, the solutions for several values of $s$. To this end, in {Figure \ref{fig:FH111}} we illustrate the solution dynamics of the fractional heat equation system \eqref{fractionalheat5} for {the four time instants $t=0.6, 1, 1.5$, and $3$. Within each of the {four} figures there are plots of solutions given for five values of $s$ and also for the PDE case. Clearly the plots show that as $s$ becomes larger, the peak of the solutions {lowers} more quickly and the spread of the solutions {becomes} larger. For example, in the plots, {the cyan, black, and blue} plots corresponding to the smaller values of $s$ have slower descent and smaller spread than do the {magenta, green, and red plots}. Of course, {comparing these plots we see that the solutions at later times have a wider spread than the corresponding solutions at earlier times.}

\begin{figure}[htb!]
    \centerline{
     \raisebox{6mm}{(a)}\includegraphics[width=0.42\linewidth]{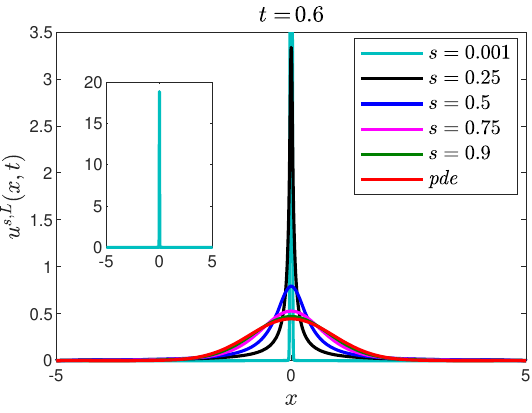}\hspace{5mm}
    \includegraphics[width=0.42\linewidth]{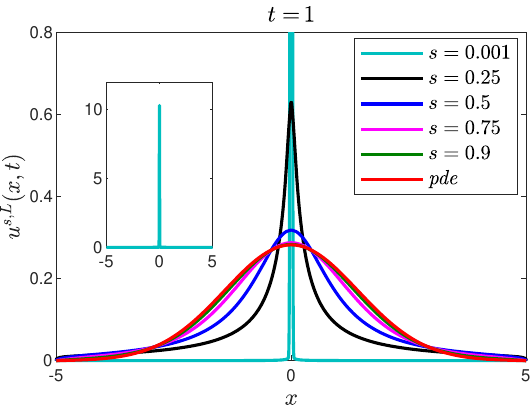}\raisebox{6mm}{\quad(b)}
  }
  \centerline{
     \raisebox{6mm}{(c)}\includegraphics[width=0.42\linewidth]{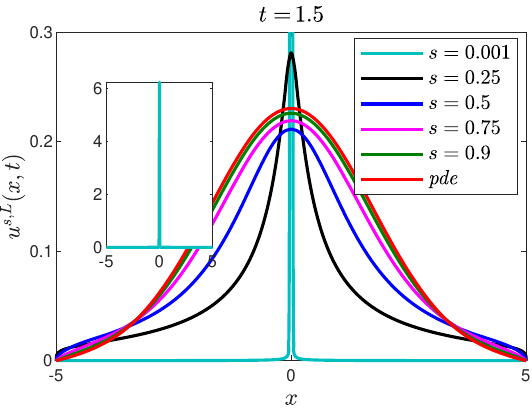}\hspace{5mm}
    \includegraphics[width=0.42\linewidth]{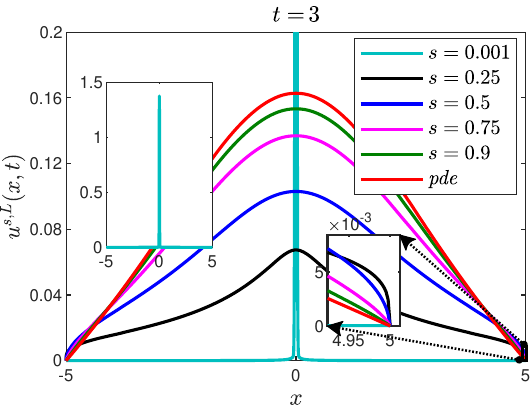}\raisebox{6mm}{\quad(d)}
  }\vspace{-3mm}
    \caption{Comparison of the solution $u^{s, L}(x, t)$ of the fractional heat system \eqref{fractionalheat5} and $u^{pde}(x, t)$ of the PDE system \eqref{classicalheat} at four time instances. To highlight the spread, only the segments within the displayed $y$-range are plotted.}
\label{fig:FH111}
\end{figure}   

\begin{re}
{\em
We end this subsection by remarking about the solution behaviors in the limits of $s\rightarrow 0^+$ and $s\rightarrow 1^-$. 
For very small values of $s$, the appreciable spread in the solution $u^{s,L}(x, t)$ becomes very narrow as is illustrated in Figure \ref{fig:fgleone}(a). Moreover, for very small $s$ we also have that the peak value of the solution $u^{s,L}(0, t)$ decays very quickly as is illustrated in Table \ref{tab:tbleone} for $s=0.001$ for six values of $t$. In this case, the fractional Laplacian operator ${\mathcal L}_s$ reduces, for all intended purposes, to the {negative of the identity operator}, as has already been noted in the first line in \eqref{fractionalheat3}. Consequently, for $s\rightarrow 0^+$, the solution of the system \eqref{fractionalheat5} can be {approximated by} 
\begin{equation}\label{szerop}
 u^{s\to 0^+,L}(x,t) = u_0(x)e^{-t}\quad\mbox{\em for any}\,\,\,\, x\in[-L,L]\,\,\,\, \mbox{\em and}\,\,\,\,t\ge0.
\end{equation}
In the limit $s\rightarrow 1^-$, the time history of the solution $u^{s\rightarrow 1^-,L}(x, t)$ approaches that of  $u^{\it pde}(x, t)$. For example, this is observed by comparing the last two entries in each of Figures \ref{tab:tbleone}(e) and \ref{tab:tbleone}(f). 
}
\qquad$\Box$
\end{re}


\section{\bf Nonlocal fractional heat equation systems}\label{sec3}

For the classical fractional equation system \eqref{fractionalheat1} given in Section \ref{sec22}, solutions are sought for all $x\in{\mathbb R}$ and domains of integration are also infinite (see \eqref{fractionalheat2}). For the system given in Section \ref{sec23}, solutions are sought on a bounded interval $x\in (-L,L)$ (see \eqref{fractionalheat5}) but the domains of integration are again infinite (see \eqref{fractionalheat4}). In this section, {\em we consider nonlocal fractional heat equation systems for which {\em both} solutions are sought in finite regions {\em and} domains of integration are also finite.}

\subsection{\bf Truncated kernel functions}
\label{sub31}

Given the kernel function $ \phi_{s}(z)$ defined in \eqref{list} we subject that function to a truncation of its spatial support and which may also be subject to the truncation of its singularity. First, given input constants $s\in(0,1)$ and $\delta>0$
we have the modified kernel function 
 \begin{equation}\label{kernel-domain-trunc}
    \phi_{s,\delta}(z) = 
    \begin{cases}
     0&\mbox{for} \,\,\, |z|> \delta
    \\[1ex]
         \displaystyle   \frac{C_{s} }{|z|^{1+2s}} \quad &\mbox{for} \,\,\, |z|\le \delta.
    \end{cases}
\end{equation}

\noindent Clearly, this kernel function has a singularity at $z=0$ and has a bounded, i.e. truncated, spatial support  $z\in[-\delta,\delta]$. 

To avoid having to deal with the singularity in \eqref{kernel-domain-trunc}, we further modify of the kernel function. Given constants $s\in(0,1)$, $\delta>0$, and $\varepsilon\le\delta$, we define the kernel function
\begin{equation}\label{kfwithgamma}
  \phi_{s,\delta,\varepsilon}(z) =
  \left\{
\begin{aligned}
&   0 && \,\,\, \mbox{for}  \,\,\, |z| > \delta
\\  
&\frac{C_s}{|z|^{1+2s}}      &&   \,\,\, \mbox{for} \,\,\,0<\varepsilon < |z| \le \delta
\\  
&\frac{C_s}{|\varepsilon|^{1+2s}}        && 
 \,\,\,\mbox{for} \,\,\,  | z|  \le  \varepsilon
\end{aligned}  \right.
\end{equation}
which is illustrated in the Figure \ref{figtrunc}.

 \begin{figure}[h]
  \centerline{  
  \includegraphics[width=0.35\linewidth, height = 5.4cm]{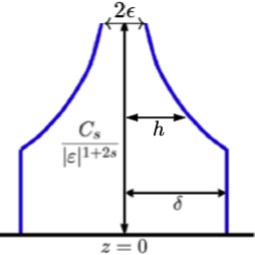}
  }
\caption{Illustration of the modified kernel function defined in \eqref{kfwithgamma}.}\label{figtrunc}
\end{figure}

It is tempting to choose {$\varepsilon=h$} or even lager values of {$\varepsilon$}. However, in our experience, we find that truncation effected with {$\varepsilon=h$} leads to poor results with respect to the accuracy of approximate solutions. Thus, as illustrated in Figure \ref{figtrunc}, we choose values of $\varepsilon\le h$. Thus the kernel function is truncated at the height ${C_s}/{|\varepsilon|^{1+2s}}$. In Table \ref{tab:C077777}  we provide the maximum value ${C_s}/{|\varepsilon|^{1+2s}}$ of the truncated kernel function $\phi_{s,\delta,\varepsilon}(z)$ for five values of $\varepsilon\le h=0.0025$ and five values of $s$. Note that the value of $s$ has a profound effect on the results given in that table, i.e.,  for a given value of $\varepsilon$ the maximum value of the kernel function increases  rapidly as the value of  $s$  increases. Also, for a given value of $s$ the  maximum value of the kernel function  decreases as the value of $\varepsilon$ increases.

\begin{table}[htb!]
\begin{center}
\begin{tabular}{|l||c|c|c|c|c|c|}
\hline
$\varepsilon =\,\longrightarrow$ & $0$ & $ 0.0001$ & $ 0.0005$  & $ 0.0010$ & $0.0020$ &$h=0.0025$ \\
\hline
$s = 0.10$ & $\infty$ & 5.698e+3   & 8.260e+2 & 3.596e+2 & 1.565e+2 & 1.197e+2 \\
$s = 0.25$ & $\infty$ & 1.995e+5   & 1.784e+4 & 6.308e+3 & 2.230e+3 & 1.596e+3 \\
$s = 0.50$ & $\infty$ & 3.183e+7   & 1.273e+6 & 3.183e+5 & 7.958e+4 & 5.093e+4\\
$s = 0.75$ & $\infty$ & 2.992e+9   & 5.353e+7 & 9.462e+6 & 1.673e+6 & 9.575e+5\\
$s = 0.90$ & $\infty$ & 2.614e+10 & 2.885e+8 & 4.142e+7 & 5.948e+6 & 3.184e+6\\
\hline
\end{tabular}
\caption{The maximal values of the kernel function { $\phi_{s, \delta, \varepsilon}(z)$} for several combinations of $s$ and  $\varepsilon$.}
\label{tab:C077777}
\end{center}
\end{table}

\subsection{\bf Nonlocal fractional heat equation systems posed on bounded integration domains}\label{sub32}

For the classical fractional diffusion system given in \eqref{fractionalheat1} not only does the operator ${\mathcal L}_{s}$ act on a function $u^{s}(x,t)$ {for all $x\in\mathbb R$} but it also {requires an integration over all $y\in\mathbb R$.}  Moreover, although for the system  \eqref{fractionalheat5} the operator { ${\mathcal L}_{s, L}$} acts on a function { $u^{s, L}(x,t)$} for $x$ in the bounded domain $(-L,L)$, it also requires an integration over all $y\in\mathbb R$. Here, we consider two alternate fractional equations systems which correspond to the two kernel functions defined in \eqref{kernel-domain-trunc} and \eqref{kfwithgamma}.

\subsubsection{\bf Kernel function having a singularity}\label{sub43}

Consider the case in which the operator ${\mathcal L}_{s,L,\delta}$ acts on a function $u^{s,L,\delta}(x,t)$ only over the bounded domain $x\in(-L,L)$ and for which we have the kernel function defined in \eqref{kernel-domain-trunc}. Thus, we have, for $s \in (0,1)$ and $t>0$, 
\begin{equation}\label{fractionalheat10}
    {\mathcal L}_{s,L,\delta} u^{s,L,\delta}(x,t) = {P.V.}\int_{x-\delta}^{x+\delta}\Big(u^{s,L,\delta}(y,t)-u^{s,L,\delta}(x,t)\Big)\phi_{s,\delta}(|y-x|)dy 
     \quad\mbox{for}\,\,\,x\in(-L,L).
\end{equation}
Then, a fractional diffusion system on bounded domains has the form
\begin{equation}\label{fractionalheat11}
\left\{
\begin{aligned}
u_t^{s,L,\delta}(x, t) = {\mathcal L}_{s,L,\delta} u^{s,L,\delta}(x, t)&\qquad\mbox{for} \ x\in  (-L,L) \,\,\,\mbox{and}\,\,\, t\in (0,T] 
\\
u^{s,L,\delta}(x,t) = g_-(x,t) &\qquad\mbox{for} \,\,\,x\in [-L-\delta,-L] \,\,\,\mbox{and}\,\,\, t\in (0,T] 
\\
u^{s,L,\delta}(x,t) = g_+(x,t) &\qquad\mbox{for} \,\,\,x\in [L,L+\delta] \,\,\,\mbox{and}\,\,\, t\in (0,T]      
\\
u^{s,L,\delta}(x,0) = u_0(x) &\qquad\mbox{for} \,\,\,x\in[-L,L]
\end{aligned}
\right.
\end{equation}
where now we impose volume constraints over {$[-L-\delta, -L] \cup [L, L+\delta]$} and the initial condition is imposed on $x\in[-L,L]$.

\subsubsection{\bf Bounded kernel function}\label{sub44}

Mimicking \eqref{fractionalheat10} and \eqref{fractionalheat11} but now using the kernel function \eqref{kfwithgamma}, we now have, for $s \in (0,1)$, $t>0$, and $x\in(-L,L)$,
\begin{equation}\label{fractionalheat1010}
    {\mathcal L}_{s,L,\delta,\varepsilon} u^{s,L,\delta,\varepsilon}(x,t) = \int_{x-\delta}^{x+\delta}\Big(u^{s,L,\delta,\varepsilon}(y,t)-u^{s,L,\delta,\varepsilon}(x,t)\Big)\phi_{s,\delta,\varepsilon}(|y-x|)dy .
\end{equation}
Then, another fractional diffusion system on bounded domains has the form
\begin{equation}\label{fractionalheat1111}
\left\{
\begin{aligned}
u_t^{s,L,\delta,\varepsilon}(x, t) = {\mathcal L}_{s,L,\delta,\varepsilon} u^{s,L,\delta,\varepsilon}(x, t)&\qquad\mbox{for} \ x\in  (-L,L) \,\,\,\mbox{and}\,\,\, t\in (0,T] 
\\
u^{s,L,\delta,\varepsilon}(x,t) = g_-(x,t) &\qquad\mbox{for} \,\,\,x\in [-L-\delta,-L] \,\,\,\mbox{and}\,\,\, t\in (0,T] 
\\
u^{s,L,\delta,\varepsilon}(x,t) = g_+(x,t) &\qquad\mbox{for} \,\,\,x\in [L,L+\delta] \,\,\,\mbox{and}\,\,\, t\in (0,T]      
\\
u^{s,L,\delta,\varepsilon}(x,0) = u_0(x) &\qquad\mbox{for} \,\,\,x\in[-L,L]
\end{aligned}
\right.
\end{equation}
where again we impose volume constraints over {$x\in[-L-\delta, -L] \cup [L, L+\delta]$} and the initial condition is imposed on $x\in[-L,L]$.

\subsection{\bf Two special cases for the systems \eqref{fractionalheat11} and \eqref{fractionalheat1111}}\label{sub43}

For each of the systems  \eqref{fractionalheat11} and \eqref{fractionalheat1111}, we  consider specific cases engendered by the two relations between $\delta$ and $L$ such that $\delta \ge 2L$  and $\delta \le L$. 

\hangone If $\delta \ge 2L$ $\Longrightarrow$ \\$\delta$ is large enough so that for every $x\in[-L,L]$  the interval $[x-L,x+L]$ is contained within the interval $[-\delta,\delta]$; this guarantees that every point in $[-L,L]$ interacts with every other point in $[-L,L]$.  Note that for this case the choice $L=\infty$ is allowable.

\noindent This case is illustrated in Figure \ref{fig:FH5555}. There, the three points $\{x_1,x_2,x_3\}\in[-L,L]$ and, in fact, for all $x\in [-L,L]$, the support of the kernel function includes all points in the red interval $[-L,L]$. Moreover, for $x\in[-L,L]$, the union of the supports of the kernel functions is $[-L-\delta,L+\delta]$ as depicted by the red, magenta, and green intervals. 
\begin{figure}[htb!]
     \centerline{
    \includegraphics[width=0.95\linewidth]{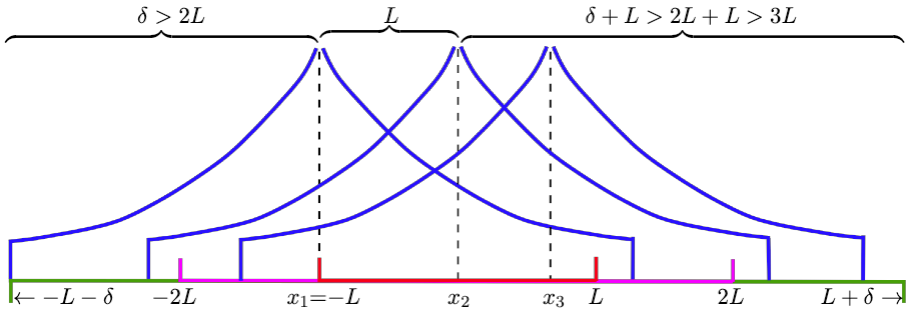}
    }
    \caption{
 For each of the three illustrations, the support of the kernel function has width $2\delta$. For each of the three specific points $x$ and for $\delta>2L$, the associated plots for each of the supports of the kernel functions is sufficiently large so that the interval $[-L,L]\subset [x-\delta,x+\delta]$.}
\label{fig:FH5555}
\end{figure}

\hangone If $\delta \le L$ $\Longrightarrow$ \\
 $\delta$ is small enough so that for every $x\in[-L,L]$ the interval $[x-\delta,x+\delta]$ is contained within the interval $[-L-\delta,L+\delta]$;  in this case, points $x\in(-L,L)$ only interact with a subset of the points $y\in[-L,L]$.
 
 \noindent This case is illustrated in Figure \ref{fig:FH6666}. There, the three points $\{x_1,x_2,x_3\}\in[-L,L]$ and, in fact, for all $x\in [-L,L]$, the support of the kernel function is too narrow to include all points in $[-L,L]$. However, for $x\in[-L,L]$, the union of the supports of the kernel functions is still $[-L-\delta,L+\delta]$ as depicted by the red and magenta intervals. 

 \begin{figure}[htb!]
     \centerline{
    \includegraphics[width=0.6\linewidth, height=5cm]{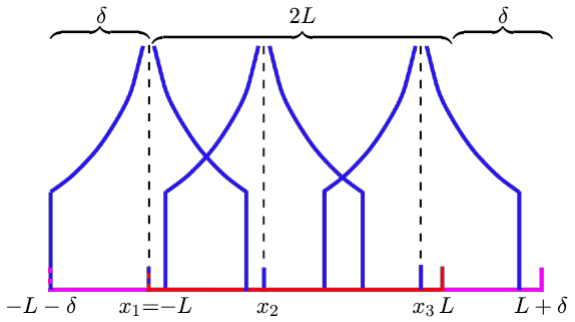}
    }
    \caption{
  For each of the three illustrations, the support of the kernel function has width $2\delta$. For each of the three specific points $x$ and for $\delta<L$, the associated plots for each of the supports of the kernel functions are sufficiently small so that the intervals $[x-\delta,x+\delta]\cap[-L,L]\not\subset [-L,L]$.}
\label{fig:FH6666}
\end{figure}

\section{Choices made about the spatial support and singularity of the kernel function}
\label{sec4}

In this section we consider four fractional diffusion systems for all of which solutions are sought in the finite domain {$\Omega=(-L,L)$} and which feature choices related to the spatial support and singularity of kernel functions. Specifically, we consider these four choices:

\hangthree-- the spatial support $[-\delta,\delta]$ of the truncated  kernel function is chosen to be { $\delta\ge 2L$}; the singularity of the kernel function is not truncated

\hangthree-- the spatial support $[-\delta,\delta]$ of the truncated  kernel function is chosen to be {$\delta \le L$}; the singularity of the kernel function is not truncated

\hangthree-- the spatial support $[-\delta,\delta]$ of the truncated  kernel function is chosen to be {$\delta\ge 2L$}; the singularity of the kernel function is truncated

\hangthree-- the spatial support $[-\delta,\delta]$ of the truncated  kernel function is chosen to be {$\delta \le L$}; the singularity of the kernel function is truncated.

\noindent The first and third of these choices (respectively the second and fourth choices) allow us to, among other information, study the limiting behavior of solutions as $\delta\to\infty$ (respectively as $\delta\to 0$).

To economize notation, we relabel these four systems according to the following table. 

\begin{table}[htb!]
\begin{center}
\begin{tabular}{|cclcll|}
\hline
superscript 
& 
system
&&
\multicolumn{1}{c}{superscript}
&
\multicolumn{1}{c}{solution}
&
$\delta,L$ relation
\\
\hline
$\{s,L,\delta\}$ 
&
\eqref{fractionalheat11}
& 
$\Longrightarrow$ 
& 
$A$ 
& $u^{A}(x,t)= {u^{s,L,\delta}(x,t)}$
& with {$\delta \ge 2L$}
\\
$\{s,L,\delta\}$
&
\eqref{fractionalheat11} 
& 
$\Longrightarrow$ 
& 
$B$
& $u^{B}(x,t)={u^{s,L,\delta}(x,t)}$
& with {$\delta \le L$}
\\
$\{s,L,\delta,\varepsilon\}$
&
\eqref{fractionalheat1111} 
& 
$\Longrightarrow$ 
& 
$C$
& $u^{C}(x,t)=u^{s,L,\delta,\varepsilon}(x,t)$
& with {$\delta \ge 2L$}
\\
$\{s,L,\delta,\varepsilon\}$
&
$\eqref{fractionalheat1111}$
& 
$\Longrightarrow$
& 
$D$
& $u^{D}(x,t)=u^{s,L,\delta,\varepsilon}(x,t)$
& with {$\delta \le L$}
\\
\hline
\end{tabular}
\end{center}
 \caption{
For sake of simplifying the exposition, this table defines the notational changes used in the rest of the paper. }
\label{fig:iiiiii}
\end{table}

\begin{re}
{\em 
Henceforth when we refer to a {\em solution} of any of these systems  we consider, we actually mean  {\em approximate solution} { obtained through spatial and temporal discretizations.}
}\quad$\Box$
\end{re}

\subsection{\bf System $A = \{s,L,\delta\}$ with $\delta\ge 2L$} \label{sys42}

In this subsection we consider the fractional heat equation system \eqref{fractionalheat11} for which the {\em kernel function has \em large \em bounded support $[-\delta,\delta]$ with $\delta\ge 2L$ but which has a \em{singularity}.} Specifically it is devoted to illustrating the convergence, as $\delta$ increases, of the solution $u^{A}(x, t)$  of the system \eqref{fractionalheat11}  to the solution $u^{s,L}(x, t)$ of the fractional heat system \eqref{fractionalheat5}. We also confirm that $u^{A}(x, t)$ does indeed model anomalous diffusion. 

Figure \ref{fig:A0} illustrates the time histories of the solution $u^{A}(x, t)$ of the system \eqref{fractionalheat11} for four values of $s$. To the naked eye, the four figures in Figure \ref{fig:A0} and the middle four figures in Figure \ref{fig:fgleone} are indistinguishable.  

\begin{figure}[p]
 \centerline{
    \raisebox{6mm}{(a)}\includegraphics[width=0.39\linewidth]{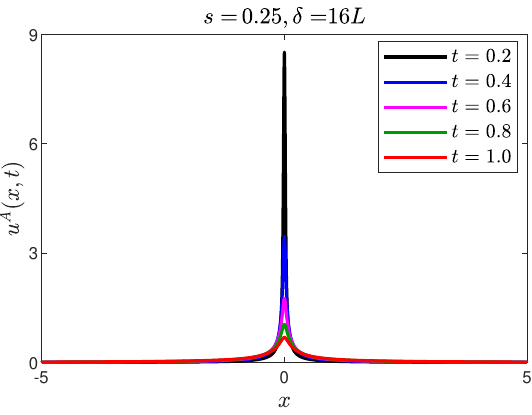}\hspace{5mm}
    \includegraphics[width=0.39\linewidth]{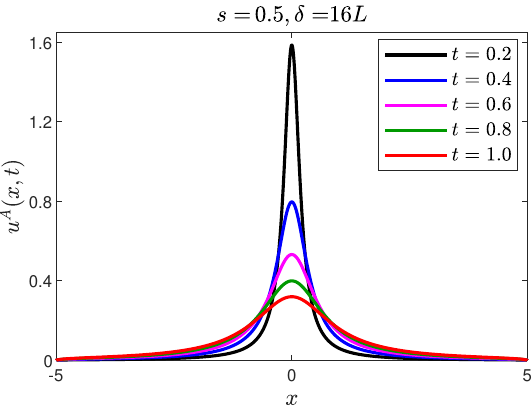}\raisebox{6mm}{\quad(b)}}
    \vspace{-1mm}
    \centerline{
    \raisebox{6mm}{(c)}
    \includegraphics[width=0.39\linewidth]{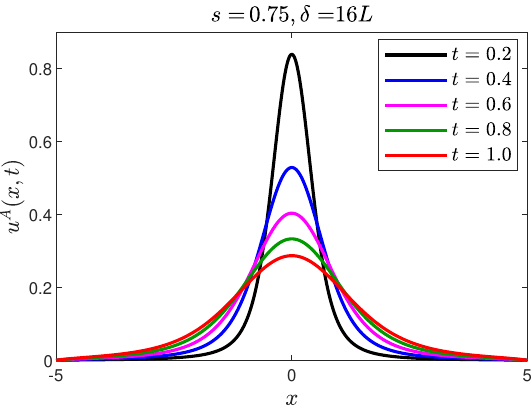}\hspace{5mm}
    \includegraphics[width=0.39\linewidth]{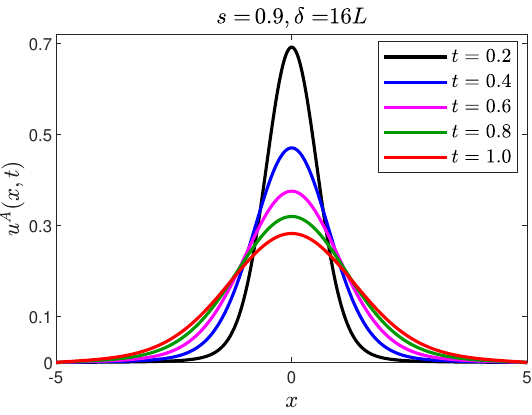}\raisebox{6mm}{\quad(d)}
    }\vspace{-3mm}
    \caption{System A: Solution dynamics with $\delta=16L$.  The four plots provide, for four values of $s$, a time history of the solution $u^{A}(x, t)$ of System A  given in \eqref{fractionalheat11}.  Note that the scale of the ordinates differ for each plot.}\label{fig:A0}
\end{figure}

\begin{figure}[p]
    \centerline{
    \raisebox{6mm}{(a)}\includegraphics[width=0.391\linewidth]{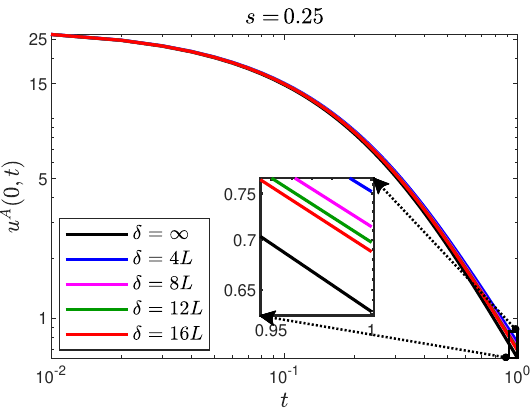}\hspace{5mm}
    \includegraphics[width=0.391\linewidth]{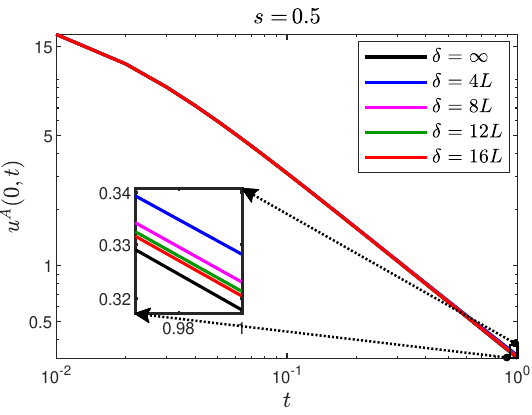}\raisebox{6mm}{\quad(b)}
    }\vspace{-1mm}
    \centerline{
    \raisebox{6mm}{(c)}\includegraphics[width=0.391\linewidth]{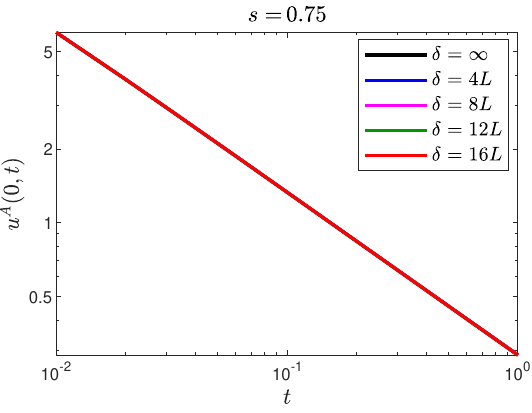}\hspace{2mm}
    \includegraphics[width=0.391\linewidth]{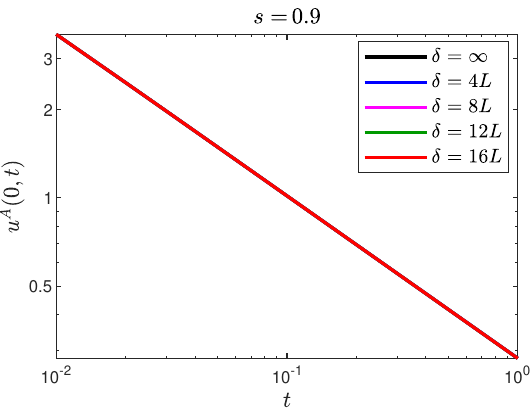}\raisebox{6mm}{\quad(d)}
    }\vspace{-3.5mm}
    \caption{System A: The time evolution of $u^{A}(0, t)$ for four values of $s$, four values of $\delta$, and $L=5$. For comparison purposes we also plot the black curves that correspond to the time evolution of  $u^{s,L}(0, t)$ for which $\delta =\infty$. 
    }\label{fig:A1}
\end{figure}

Figures \ref{fig:A1} and  \ref{fig:A3} respectively illustrate, for four values of $s$, the convergence, as the value of $\delta$ increases, of the solution $u^{A}(0, t)\to u^{s,L}(0, t)$ and of the $\ell^2$-norm of the solution $\|u^{A}(\cdot, t)\|_{\ell_2} \to \|u^{s,L}(\cdot, t)\|_{\ell_2}$. The  inserts in Figures \ref{fig:A1} and  \ref{fig:A3} illustrate that the convergence behavior is slower for smaller values of $s$ compared to larger values of $s$. Overall, a perusal of these two figures and of Figure \ref{fig:A0} indicate that  System A provides a good approximation to the solution of fractional heat system \eqref{fractionalheat5} when $\delta$ is large enough.

\begin{figure}[h!]
    \centerline{ \raisebox{6mm}{(a)}\includegraphics[width=0.391\linewidth]{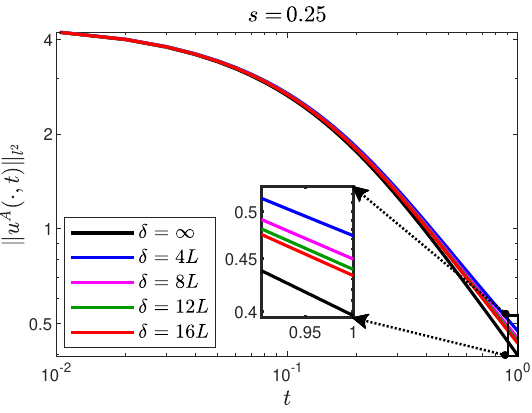}\hspace{5mm}
    \includegraphics[width=0.391\linewidth]{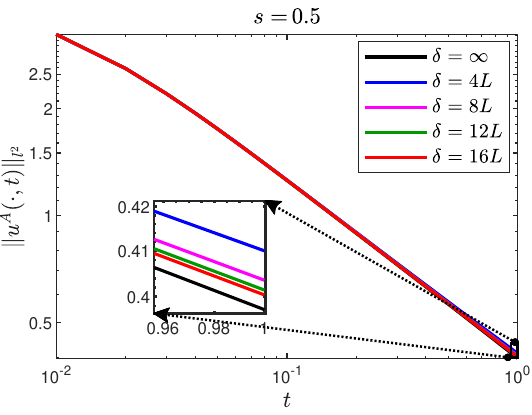}\raisebox{6mm}{\quad(b)}
    }\vspace{-1mm}
    \centerline{
    \raisebox{6mm}{(c)}\includegraphics[width=0.391\linewidth]{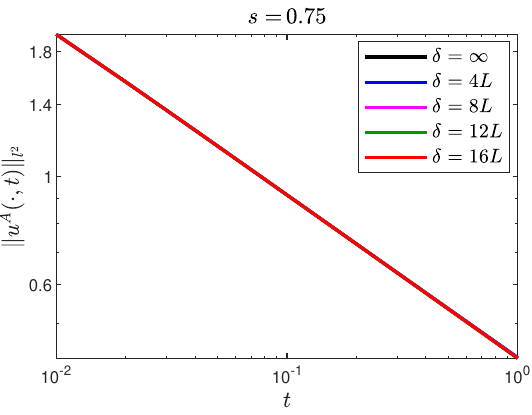}\hspace{5mm}
    \includegraphics[width=0.391\linewidth]{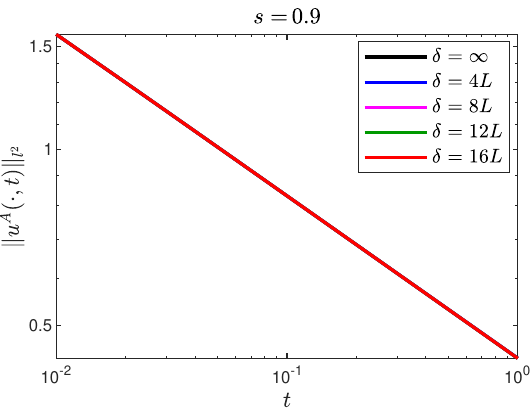}
    \raisebox{6mm}{\quad(d)}
    }\vspace{-3.5mm}
\caption{System A: The time evolution of $\|u^{A}(\cdot, t)\|_{\ell_2}$ for four values of $s$, four values of $\delta$, and $L=5$. For comparison purposes we also plot the black curves that correspond to the evolution of  $\|u^{s,L}(\cdot,t)\|_{\ell^2}$ for which $\delta =\infty$. 
  }\label{fig:A3}
\end{figure}

The observations just made are further supported by Table \ref{tab:A4} which provides, for four values of $s$, the error $\|u^{A}(\cdot, t)-u^{s,L}(\cdot, t)\|_{\ell_2}$ at $t = 1$ and the corresponding convergence rates with respect to $\delta$.\footnote{Here, the convergence rate is computed by $\log\big(\|u^{A_1}(\cdot, t)-u^{s,L}(\cdot, t)\|_{\ell^2}/\|u^{A_2}(\cdot, t)-u^{s,L}(\cdot, t)\|_{\ell^2}\big)/\log(\delta_1/\delta_2)$ for two consecutive $\delta_1$ and $\delta_2$, where $A_1=\{s,L,\delta_1\}$ and $A_2=\{s,L,\delta_2\}$.}   The table further illustrates that the solution of System A converges, as $\delta$ increases, to that of the fractional heat equation \eqref{fractionalheat5}, which is consistent with our observations about Figures \ref{fig:A1} and \ref{fig:A3}. Moreover, note that the theoretical decay rate of the error (as  proven in \cite{Burkovska2019}) is given as
\begin{equation}\label{sexp}
\|u^{A}(\cdot, t)-u^{s,L}(\cdot, t)\|_{\ell_2} \sim {\mathcal O}(\delta^{-2s}) \qquad \mbox{as \ $\delta \to \infty$}
\end{equation}
so that the approximate decay rates given in Table  \ref{tab:A4} are in very close agreement with the theoretical rates given in  \eqref{sexp}. An interesting observation is that to achieve the same error, smaller values of $s$ require a larger $\delta$ compared to that for a larger value of $s$.
In Table  \ref{tab:A4}, for each value of $s$ the {\em error decreases monotonically} as $\delta$ increases and for each value of $\delta$ the {\em error decreases monotonically} as $s$ increases.

\begin{table}[htb!]
    \centering
    \begin{tabular}{@{}|l||cccccccc|@{}} 
    \hline
     &\multicolumn{2}{c}{$s = 0.25$}  &\multicolumn{2}{c}{$s = 0.50$}
     &\multicolumn{2}{c}{$s = 0.75$} &\multicolumn{2}{c|}{$s = 0.90$} \\ 
    \cline{2-9}
    & error & rate & error & rate & error & rate & error & rate  \\
    \hline 
$\delta = 8L$  &5.328e-2 & -- -- &6.372e-3 &-- --  &6.708e-4 & -- -- &1.052e-4 &-- --\\
$\delta = 10L$ &4.733e-2 &-0.530 &5.090e-3 &-1.007 &4.799e-4 &-1.501 &7.040e-5 &-1.798\\
$\delta = 12L$ &4.299e-2 &-0.528 &4.237e-3 &-1.006 &3.651e-4 &-1.500 &5.073e-5 &-1.797\\
$\delta = 14L$ &3.965e-2 &-0.525 &3.629e-3 &-1.005 &2.897e-4 &-1.500 &3.847e-5 &-1.796\\
$\delta = 16L$ &3.697e-2 &-0.523 &3.174e-3 &-1.004 &2.371e-4 &-1.500 &3.027e-5 &-1.793\\
$\delta = 32L$ &2.580e-2 &-0.519 &1.584e-3 &-1.003 &8.391e-5 &-1.499 &8.858e-6 &-1.773\\
    \hline
anal. rate & & -0.5  & & -1.0 & & -1.5 & &-1.8\\
    \hline
\end{tabular}
\caption{System A: For $L=5$, the numerically obtained approximate errors $\|u^{A}(\cdot, t)-u^{s,L}(\cdot, t)\|_{\ell^2}$  and the convergence rates with respect to $\delta$. The last row provides the theoretical decay rates in  \eqref{sexp}.}\label{tab:A4}
\end{table}


\subsection{\bf System $B= \{s,L,\delta\}$ with $\delta\le L$}\label{sys43}

In this subsection we consider the system \eqref{fractionalheat11} for which the {\em kernel function has \em small \em support $[-\delta,\delta]$ with $\delta\le L$ but which has a \em{singularity}.} Specifically, it is devoted to illustrating the relations between the solution $u^{B}(x, t)$  of the system \eqref{fractionalheat11} to the solution $u^{pde}(x, t)$ of the partial differential equation system \eqref{classicalheat} and also to the solution $u^{s,L}(x, t)$ of the fractional system \eqref{fractionalheat5}. 

Figure \ref{fig:A2} illustrates the time histories of the solution $u^{B}(x, t)$ of the system \eqref{fractionalheat11} for fixed $L=5$ and $\delta = L/8 = 0.625$, and for four values of $s$. We note the visually excellent similarity between Figure \ref{fig:A2}(d) with Figure \ref{fig:fgleone}(f). Naturally, we expect such a similarity because $s=0.98$ in Figure \ref{fig:A2}(d). 

\begin{figure}[htb!]
    \centerline{
  \raisebox{6mm}{(a)}
  \includegraphics[width=0.391\linewidth]{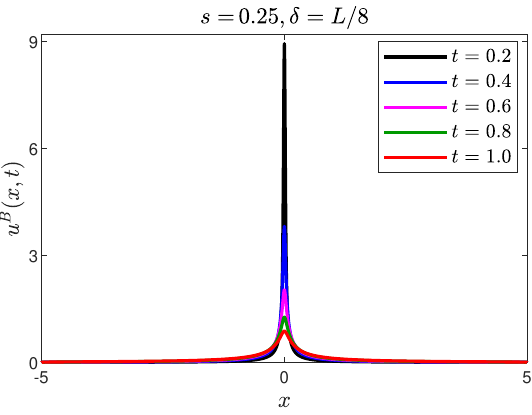}\hspace{5mm}
  \includegraphics[width=0.391\linewidth]{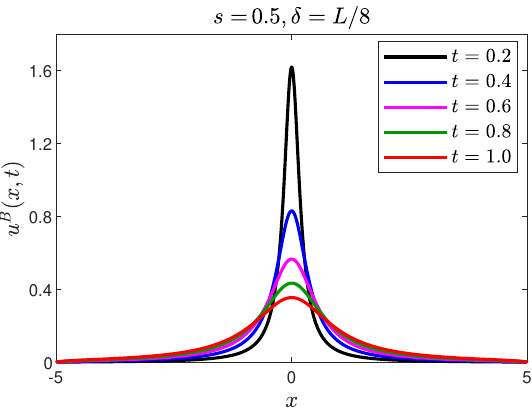}\raisebox{6mm}{\quad(b)}
    }\vspace{-1mm}
    \centerline{
  \raisebox{6mm}{(c)}
  \includegraphics[width=0.391\linewidth]{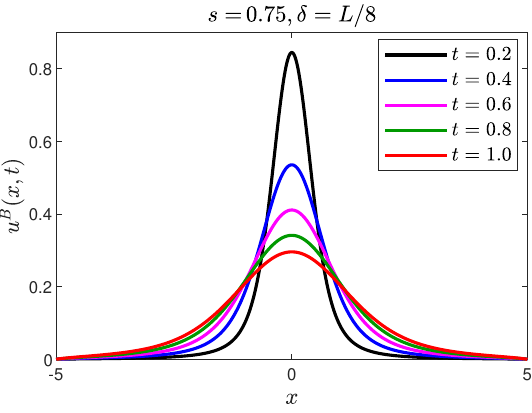}\hspace{5mm}
  \includegraphics[width=0.391\linewidth]{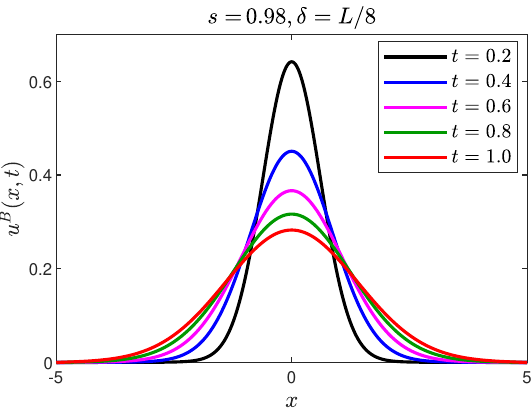}\raisebox{6mm}{\quad(d)}
    }\vspace{-3mm}
    \caption{System B: Solution dynamics with $\delta=L/8$.  The four plots provide, for four values of $s$, a time history of the solution {$u^{B}(x, t)$} of System B  given in \eqref{fractionalheat11}. 
    When examining these plots, one should include the first and last plots in Figure \ref{fig:fgleone}.}\label{fig:A2}
\end{figure}

Table \ref{tab:B2} shows, for small values of $\delta$, the differences between the solution  $u^{B}(x, t)$ of System B and the solution  $u^{pde}(x, t)$ of classical heat equation.  Thus the convergence  $u^{B}(x, t)$ to $u^{pde}(x, t)$ as $s\to1^-$ is numerically verified in that table. Note that as $s$ increases the error decreases whereas for a fixed value of $s$ and for very small $\delta$, the error remains virtually constant. 
\begin{table}[htb!]
\begin{center}
\begin{tabular}{|l||cccccc|} \hline  
            $\delta =\, \longrightarrow$ &$0.6250$ 
            &$0.0400$  &$0.0200$ &$0.0100$ &$0.0050$  &$0.0025$\\ \hline
\multirow{1}{*}{$s=0.75$}   &3.343e-2 &3.340e-2 &3.340e-2 &3.340e-2 &3.340e-2  & 3.340e-2 \\
\multirow{1}{*}{$s=0.90$}  &1.283e-2 &1.153e-2 &1.151e-2 &1.151e-2 &1.151e-2 &1.151e-2  \\
\multirow{1}{*}{$s=0.99$}  &1.194e-3  &1.068e-3 &1.067e-3 &1.066e-3 &1.066e-3 &1.066e-3  \\
\multirow{1}{*}{$s=0.999$} &1.187e-4 &1.062e-4 &1.060e-4 &1.060e-4 &1.060e-4 &1.059e-4  \\\hline
\end{tabular}
\caption{System B: For $L=5$ and for several $\delta$ and for $s\to1^-$, the numerical errors $\|u^B(\cdot, t)-u^{pde}(\cdot, t)\|_{\ell_2}$ at time $t  = 1$.}
\label{tab:B2}
\end{center}
\end{table}


\subsection{\bf System $C = \{s,L,\delta,\varepsilon\}$ \bf with $\delta\ge 2L$ and $\varepsilon\ll 1$}\label{sys44}

In this subsection we consider the fractional heat equation system \eqref{fractionalheat11} for which the {\em kernel function has \em large \em bounded support $[-\delta,\delta]$ with $\delta\ge 2L$ and for which the singularity of the kernel function is truncated.}
Specifically it is devoted to illustrating the convergence, as $\delta$ increases and $\varepsilon$ decreases, of the solution $u^{C}(x, t)$  of the system \eqref{fractionalheat11}  to the solution $u^{s,L}(x, t)$ of the fractional heat system \eqref{fractionalheat5}. We also confirm that $u^{C}(x, t)$ does indeed model anomalous diffusion. 

Figure \ref{fig:CC2} provides plots of the time histories of $u^{C}(x, t)$ for four values of $s$,  $\varepsilon=h=0.0025$, $L=5$, and $\delta=64L$ so that $\delta=320$.  Note that the choice of $\varepsilon$ is independent of the numerical parameter $h$, although in this figure we present the results with $\varepsilon = h$.
One observes that the temporal decays of the solutions become greater as the value of $s$ increases.

\begin{figure}[htb!]
    \centerline{
    \raisebox{6mm}{(a)}
    \includegraphics[width=0.391\linewidth]{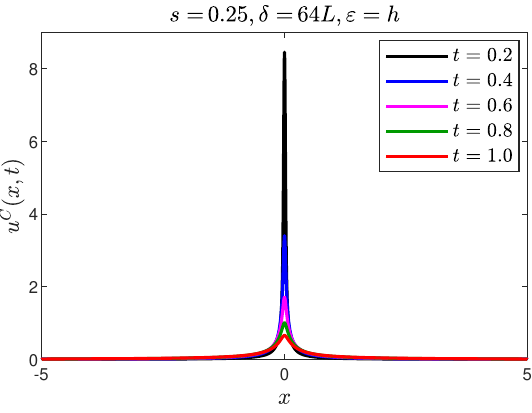}\hspace{5mm}
    \includegraphics[width=0.391\linewidth]{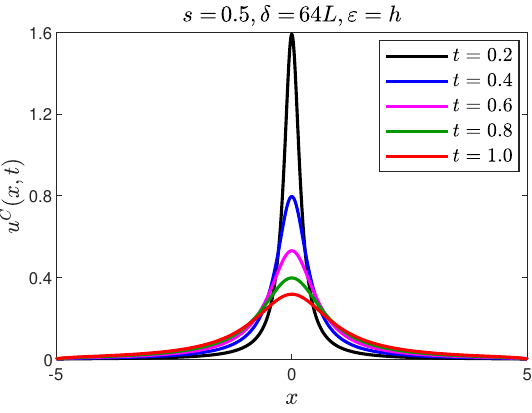}\raisebox{6mm}{\quad(b)}
    }\vspace{-1mm}
    \centerline{
    \raisebox{6mm}{(c)}\includegraphics[width=0.391\linewidth]{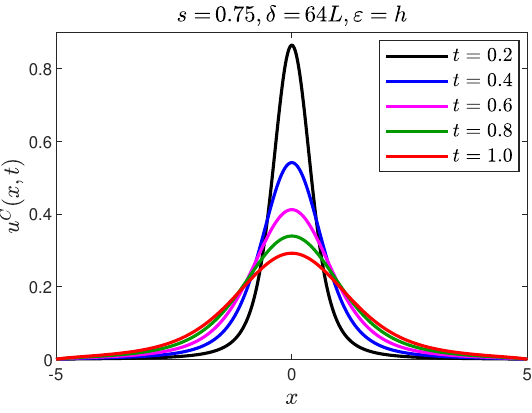}\hspace{5mm}
    \includegraphics[width=0.391\linewidth]{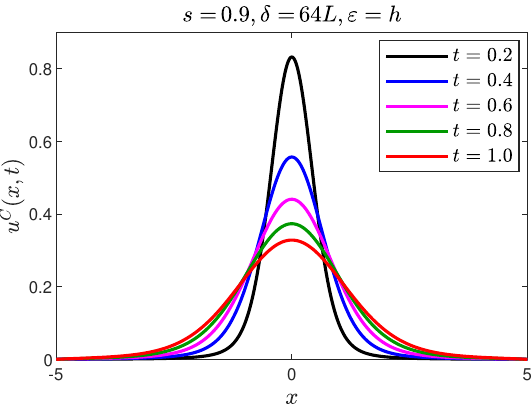}\raisebox{6mm}{\quad(d)}
    }\vspace{-3.5mm}
    \caption{System C: Solution dynamics with $\delta=64L$.  The four plots provide, for four values of $s$ and for $\varepsilon=h$, a time history of the solution{$u^{C}(x, t)$} of System C  given in \eqref{fractionalheat11}. When examining these plots, one should include the first and last plots in Figure \ref{fig:fgleone}. Note that the scale of the ordinates differ for each plot.}
 \label{fig:CC2}
\end{figure}

It is interesting to compare the solutions provided in Figures \ref{fig:A0} and \ref{fig:CC2} because one sees that the introduction of the truncation of the singularity of the kernel function results in a difference between the two solutions. For example, one notices that the decays of the peaks and the extents of the spreads of the two solutions are different. However, both Systems A and C result in anomalous diffusion.

Figures \ref{fig:C1} and \ref{fig:C3} respectively provide plots the temporal evolution of $u^C(0,t)$  and $\|u^C(\cdot,t)\|_{\ell^2}$ for four values of $s$ and five values of $\varepsilon$ and for $\delta=64L=320$ and $h=0.0025$. A visual look at these three figures give credence that $u^{C}(x, t)$ does indeed model anomalous diffusion. 
They further demonstrate that, in approximating the fractional heat equation \eqref{fractionalheat5}, both $\delta$ and $\varepsilon$ should be chosen depending on $s$. For instance, with a fixed $\varepsilon = h$ and time $t$, Figure~\ref{fig:C1} shows that the difference $|u^C(0, t) - u^{s, L}(0, t)|$ increases as $s$ increases, indicating that a smaller $\varepsilon$ is required for larger values of $s$.

\begin{figure}[htb!]
    \centerline{
    \raisebox{6mm}{(a)}
    \includegraphics[width=0.391\linewidth]{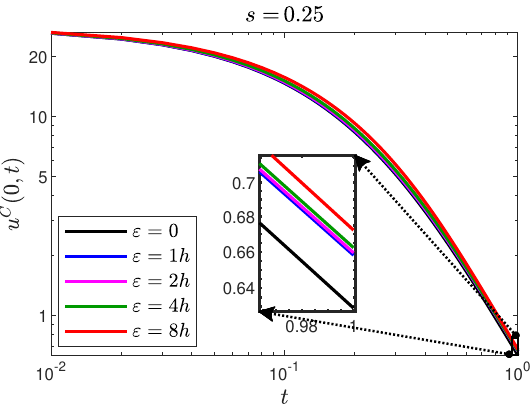}\hspace{5mm}
    \includegraphics[width=0.391\linewidth]{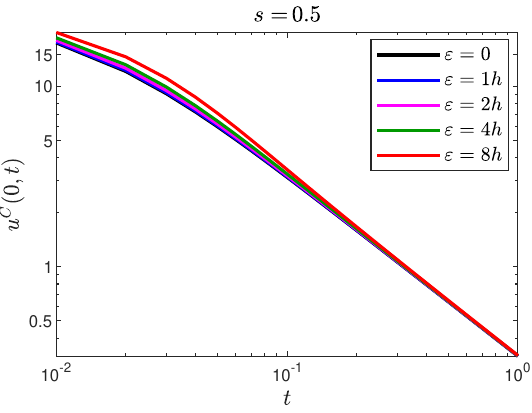}
    \raisebox{6mm}{\quad(b)}
    }\vspace{-1mm}
    \centerline{
    \raisebox{6mm}{(c)}\includegraphics[width=0.391\linewidth]{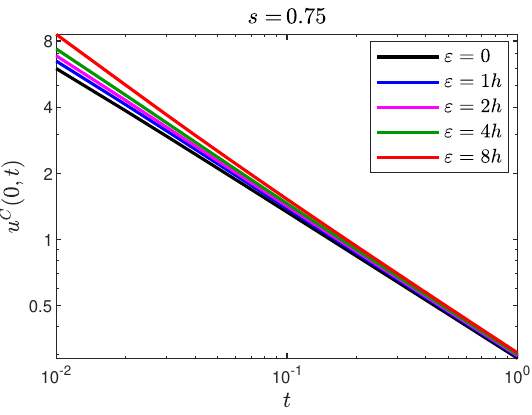}\hspace{5mm}
    \includegraphics[width=0.391\linewidth]{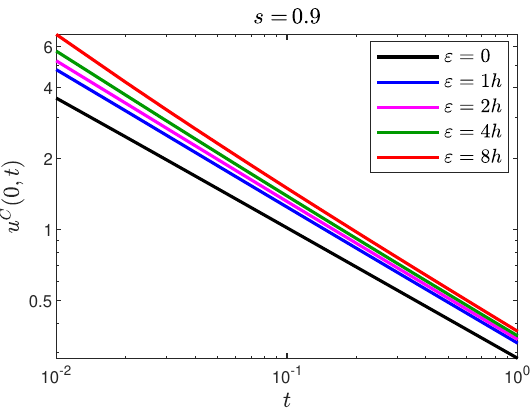}\raisebox{6mm}{\quad(d)}
    }\vspace{-3.5mm}
    \caption{System C: The time evolution of{$u^{C}(0, t)$} for four values of $s$ and for four values of $\varepsilon$ that are multiples of $h$ with $h=0.0025$ and $\delta=64L$ with $L=5$. For comparison purposes we also plot the black curves that correspond to the time evolution of  $u^{s,L}(0,t)$ for which $\varepsilon =0$. Note that the scale of the ordinates differ for each plot.}\label{fig:C1}
\end{figure}

\begin{figure}[htb!]
    \centerline{
    \raisebox{6mm}{(a)}
    \includegraphics[width=0.391\linewidth]{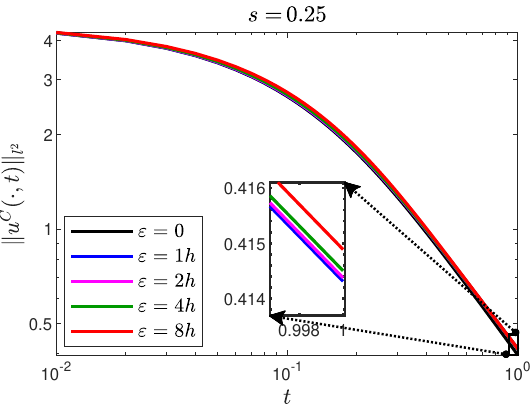}\hspace{5mm}
    \includegraphics[width=0.391\linewidth]{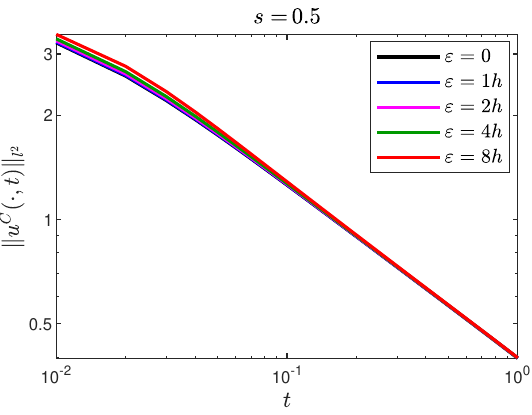}\raisebox{6mm}{\quad(b)}
    }\vspace{-1mm}
    \centerline{
    \raisebox{6mm}{(c)}\includegraphics[width=0.391\linewidth]{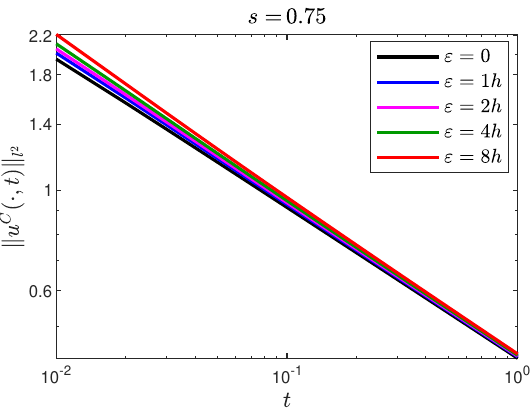}\hspace{5mm}
    \includegraphics[width=0.391\linewidth]{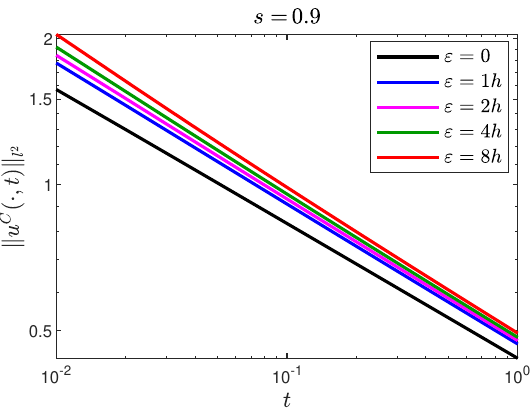}\raisebox{6mm}{\quad(d)}
    }\vspace{-3mm}
    \caption{System C: The time evolution of {$\|u^{C}(\cdot, t)\|_{\ell^2}$} for four values of $s$ and for four values of $\varepsilon$ that are multiples of $h$ with $h=0.0025$ and $\delta=64L$ with $L=5$. For comparison purposes we also plot the black curves that correspond to the time evolution of  $\|u^{s,L}(\cdot,t)\|_{\ell^2}$ for which $\varepsilon =0$. Note that the scale of the ordinates differ for each plot. }\label{fig:C3}
\end{figure}

Our results above demonstrate that System~C indeed models anomalous diffusion, 
with $\varepsilon$ and $\delta$ serving as two model parameters that depend on the underlying physical problem. 
In addition, System~C can be used to approximate the fractional heat equation~\eqref{fractionalheat5}. 
In this context, the error between System~C and the exact solution of \eqref{fractionalheat5} consists of two components: truncation errors of the  support domain and truncation errors of the singularity.

For a prescribed tolerance $\mathrm{tol}$, the parameters of $\delta$ and $\varepsilon$ should satisfy
\begin{equation}\label{eq:tol_constraints}
    \delta \leq \delta_{\mathrm{tol}}^s := c_1 \, \mathrm{tol}^{-\frac{1}{2s}}, 
    \qquad  
    \varepsilon \leq \varepsilon_{\mathrm{tol}}^s := c_2 \, \mathrm{tol}^{\tfrac{1}{2(1-s)}},
\end{equation}
with constants $c_1$ and $c_2$ that may depend on $s$. 
These conditions provide practical guidance for selecting appropriate values of $\delta$ and $\varepsilon$ 
so that System~C yields an accurate approximation to the fractional heat equation~\eqref{fractionalheat5}.

\vskip10pt

\noindent\underline{\em Convergence with respect  to $\delta$ with {$\varepsilon=\varepsilon^s_{\rm tol}$}}

Table \ref{tab:CCCC41c} provides, for four values of $s$, the error $\|u^{C}(\cdot, t)-u^{s,L}(\cdot, t)\|_{\ell_2}$ at $t = 1$ and the corresponding convergence rates with respect to $\delta$.  
Note that the errors come from {the  truncation of both} support domain and the singularity. 
Here, the singularity is truncated with {$\varepsilon = \varepsilon_{\mathrm{tol}}^s \ll h$}, ensuring that the errors from singularity truncation are sufficiently small and remain separate from those due to the truncation of support domain.  
For easy comparison, the last row provides the errors of $\delta = \infty$ and $\varepsilon = \varepsilon_{\mathrm{tol}}^s$, i.e. only the truncation of the singularity is considered.
The table shows that for each $s$, the solution of System C converges to that of the fractional heat system \eqref{fractionalheat5} as $\delta$ increases. 
Moreover, we numerically observe that the decay rate of the error is given as
\begin{equation}\label{sexp2}
\|u^{C}(\cdot, t)-u^{s,L}(\cdot, t)\|_{\ell_2} \sim {\mathcal O}(\delta^{-2s}) \qquad \mbox{as \ $\delta \to \infty$}
\end{equation}
so that the decay rates given in Table \ref{tab:CCCC41c} are in very close agreement with the expected rates given in \eqref{sexp}. 

\begin{table}[htb!]
\begin{center}
    \centering
    \begin{tabular}{@{}|l||cccccccc|@{}} 
    \hline
     &\multicolumn{2}{c}{$s = 0.25$}  &\multicolumn{2}{c}{$s = 0.50$}
     &\multicolumn{2}{c}{$s = 0.75$} &\multicolumn{2}{c|}{$s = 0.90$} \\ 
    \cline{2-9}
    & error & rate & error & rate & error & rate & error & rate  \\
    \hline
$\delta = 8L$  &5.331e-2 & -- -- &6.428e-3 &-- --  &6.715e-4 & -- -- &1.052e-4 &-- --\\
$\delta = 10L$ &4.736e-2 &-0.530 &5.146e-3 &-0.997 &4.806e-4 &-1.499 &7.049e-5 &-1.796\\
$\delta = 12L$ &4.302e-2 &-0.527 &4.293e-3 &-0.994 &3.657e-4 &-1.498 &5.082e-5 &-1.794\\
$\delta = 14L$ &3.968e-2 &-0.525 &3.685e-3 &-0.990 &2.904e-4 &-1.497 &3.856e-5 &-1.792\\
$\delta = 16L$ &3.700e-2 &-0.523 &3.230e-3 &-0.987 &2.378e-4 &-1.496 &3.037e-5 &-1.788\\
    \hline
$\delta = \infty$ &1.343e-4  &-- --&1.305e-4 & -- -- &2.719e-6 & -- -- &1.537e-6 & -- -- \\ 
\hline
\end{tabular}
\caption{System C: For $L=5$, {$\varepsilon=\varepsilon_{\rm tol}^s$}, and for four values of $s$, the numerically obtained approximate errors $\|u^{C}(\cdot, t)-u^{s,L}(\cdot, t)\|_{\ell^2}$ at time $t  = 1$ and the convergence rates with respect to $\delta$. The last row provides the errors of $\delta = \infty$ and $\varepsilon = \varepsilon_{\rm tol}^s$, i.e., only the truncation of the singularity is considered.
}\label{tab:CCCC41c}
\end{center}
\end{table}

\vskip10pt

\noindent
\underline{\em Convergence with respect to $\varepsilon$ with {$\delta = \delta_{\rm tol}^s$}}

Table \ref{tab:C2} provides, for four values of $s$, the error $\|u^{C}(\cdot, t)-u^{s,L}(\cdot, t)\|_{\ell_2}$ at $t = 1$ and the corresponding convergence rates with respect to $\varepsilon$.    
Here, the support domain is truncated with {$\delta = \delta_{\rm tol}^s$}, ensuring that the errors from support domain truncation are sufficiently small and remain separate from those due to the singularity truncation. 
{We also set the grid size to 
$h = 0.00125$, so the discretization error is negligible compared with the singularity-truncation error.}
For easy comparison, the last row provides the errors of $\varepsilon = 0$ and $\delta = \delta_{\rm tol}^s$, i.e. only the truncation of the support domain is considered. 
The table shows that as $\varepsilon$ decreases, the errors are reduced, indicating that System C provides a good approximation to the fractional heat equation \eqref{fractionalheat5}. 
Moreover, our numerical observations indicate that the convergence rates well approximate the formula
\begin{equation}\label{eps2}
\|u^{C}(\cdot, t)-u^{s,L}(\cdot, t)\|_{\ell_2} \sim {\mathcal O}(\varepsilon^{2(1-s)}) \qquad \mbox{as $\varepsilon \to 0$}
\end{equation}
Hence, to achieve the same error, a smaller $\varepsilon$ is required for a lager $s$.

\begin{table}[htb!]
\begin{center}
\begin{tabular}{@{}|l||cccccccc|@{}} 
    \hline
     &\multicolumn{2}{c}{$s = 0.25$}  &\multicolumn{2}{c}{$s = 0.50$}
     &\multicolumn{2}{c}{$s = 0.75$} &\multicolumn{2}{c|}{$s = 0.90$} \\ 
    \cline{2-9}
    & error & rate & error & rate & error & rate & error & rate  \\
    \hline
$\varepsilon = 1/800$  &6.116e-5 &-- -- &1.635e-4 &-- -- &4.781e-3 & -- --&5.124e-2 &-- --\\
$\varepsilon = 1/1600$ &2.099e-5 &1.543 &8.178e-5 &0.999 &3.360e-3 &0.509 &4.350e-2 &0.236\\
$\varepsilon = 1/2400$ &1.108e-5 &1.574 &5.460e-5 &0.997 &2.736e-3 &0.507 &3.961e-2 &0.232\\
$\varepsilon = 1/3200$ &6.997e-6 &1.599 &4.102e-5 &0.994 &2.366e-3 &0.505 &3.708e-2 &0.290\\
$\varepsilon = 1/4000$ &4.902e-6 &1.595 &3.287e-5 &0.992 &2.114e-3 &0.505 &3.525e-2 &0.272\\
$\varepsilon = 1/4800$ &3.697e-6 &1.545 &2.745e-5 &0.989 &1.928e-3 &0.504 &3.382e-2 &0.260\\
    \hline
$\varepsilon = 0$ &1.987e-6 & -- -- &1.675e-6 & -- -- &2.105e-6 & -- -- &8.100e-6 & -- -- \\ 
\hline
\end{tabular}
\caption{System C: For $L=5$,  $\delta=\delta_{\rm tol}^s$, and for four values of $s$, the numerically obtained approximate errors $\|u^{C}(\cdot, t)-u^{s,L}(\cdot, t)\|_{\ell^2}$ at time $t  = 1$ and the convergence rates with respect to $\varepsilon$. The last row provides the errors of $\varepsilon = 0$ and $\delta = \delta_{\rm tol}^s$, i.e., only the truncation of the support domain is considered. }\label{tab:C2}
\end{center}
\end{table}
}


\subsection{\bf System $D = \{s,L,\delta,\varepsilon\}$ \bf with $\delta\le L$ and $\varepsilon\ll 1$} \label{sys45}

Because the content of this subsection is related to Section \ref{sys44} in the same way that Section \ref{sys43} is related to Section  \ref{sys42}, here we only provide a very shortened consideration of System D. 
{In fact, we have one to one correspondences between Figure \ref{fig:D1} and Figure \ref{fig:A2} and between Table \ref{tab:D1} and  Table \ref{tab:B2}. Of course, there is a very stark difference between System D and System B, namely that for the latter system, the singularity of the kernel function is not truncated. 

\begin{figure}[htb!]
    \centerline{
    \raisebox{6mm}{(a)}
    \includegraphics[width=0.389\linewidth]{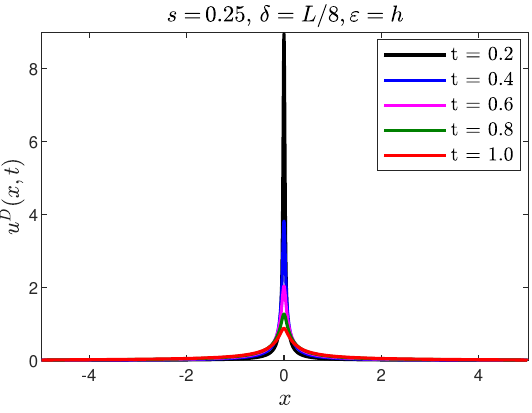}\hspace{5mm}
    \includegraphics[width=0.391\linewidth]{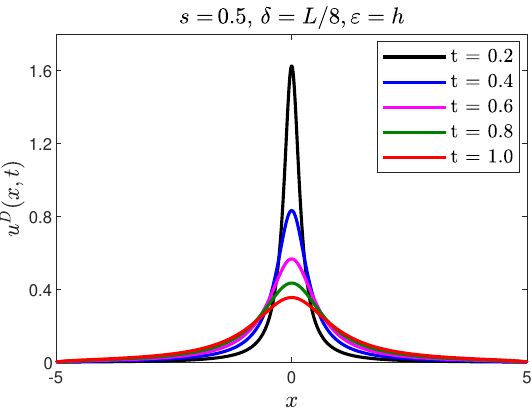}\raisebox{6mm}{\quad(b)}
    }\vspace{-1mm}
    \centerline{
    \raisebox{6mm}{(c)}\includegraphics[width=0.391\linewidth]{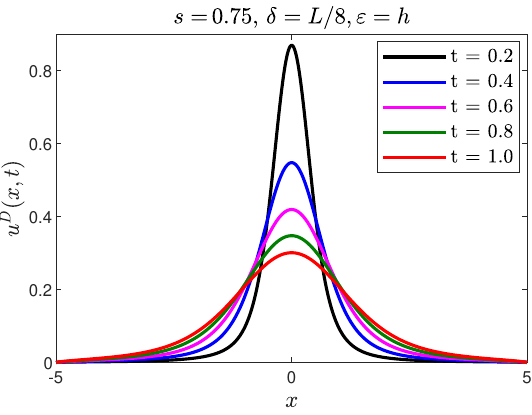}\hspace{5mm}
    \includegraphics[width=0.391\linewidth]{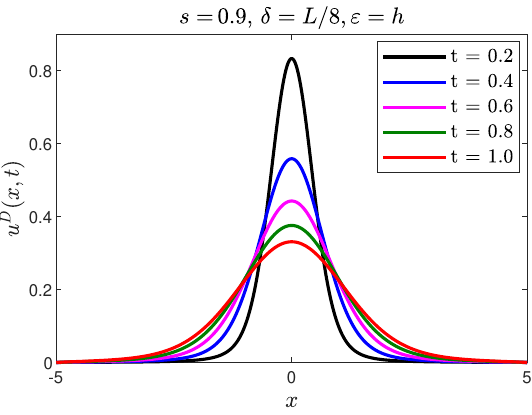}\raisebox{6mm}{(d)}
    }\vspace{-3mm}
    \caption{System D: Solution dynamics with $\delta=L/8$.  The four plots provide, for four values of $s$ and for $\varepsilon=h$, a time history of the solution $u^{D}(x, t)$ of System D. When examining these plots, one should include the first and last plots in Figure \ref{fig:fgleone}. Note that the scale of the ordinates differ for each plot.}
    \label{fig:D1}
\end{figure}

Table \ref{tab:D1} shows, for small values of $\delta$, the differences between the solution  $u^{B}(x, t)$ of System D and the solution  $u^{pde}(x, t)$ of classical heat equation, where the singularity is truncated with $\varepsilon = \varepsilon_{tol}^s$.  Thus the convergence  $u^{D}(x, t)$ to $u^{pde}(x, t)$ as $s\to1^-$ is numerically verified in that table. Note that as $s$ increases the error decreases whereas for a fixed value of $s$ and for very small $\delta$, the error remains virtually constant. 
\begin{table}[htb!]
\begin{center}
\begin{tabular}{|l||cccccc|} \hline  
            $\delta =\, \longrightarrow$ &$0.6250$ 
            &$0.0400$  &$0.0200$ &$0.0100$ &$0.0050$  &$0.0025$\\ \hline
\multirow{1}{*}{$s=0.75$}  &3.396e-2 &3.343e-2 &3.341e-2 &3.340e-2 &3.339e-2 &3.339e-2 \\
\multirow{1}{*}{$s=0.90$}  &1.181e-2 &1.153e-2 &1.151e-2 &1.151e-2 &1.151e-2 &1.151e-2  \\
\multirow{1}{*}{$s=0.99$}  &1.098e-3 &1.068e-3 &1.066e-3 &1.066e-3 &1.066e-3 &1.066e-3  \\
\multirow{1}{*}{$s=0.999$} &1.0917e-4 &1.0616e-4 &1.060e-4 &1.060e-4 &1.060e-4 &1.059e-4  \\\hline
\end{tabular}
\caption{System D: For $L=5$, {$\varepsilon = \varepsilon_{tol}^s$,} and for several $\delta$ and for $s\to1^-$, the numerical errors $\|u^D(\cdot, t)-u^{pde}(\cdot, t)\|_{\ell_2}$ at time $t  = 1$.}
\label{tab:D1}
\end{center}
\end{table}

\bigskip
\noindent{\bf  Acknowledgements.} M. Gunzburger thanks the US Department of Energy for their support under grant number DE-SC0023171. Y. Zhang thanks the US National Science Foundation for their support under grant number DMS-1953177.


\bibliographystyle{plain}
\bibliography{References}

\end{document}